\documentclass[journal,twoside,web]{ieeecolor}
\usepackage{generic}
\usepackage{cite}
\usepackage{amsmath,amssymb,amsfonts}
\usepackage{algorithmic}
\usepackage{graphicx}
\usepackage{textcomp}
\def\BibTeX{{\rm B\kern-.05em{\sc i\kern-.025em b}\kern-.08em
    T\kern-.1667em\lower.7ex\hbox{E}\kern-.125emX}}

\usepackage{savesym}
\savesymbol{AND}

\usepackage[ruled,linesnumbered]{algorithm2e}
\makeatletter
\newcommand{\nosemic}{\renewcommand{\@endalgocfline}{\relax}}
\newcommand{\dosemic}{\renewcommand{\@endalgocfline}{\algocf@endline}}
\let\oldnl\nl
\newcommand{\nonl}{\renewcommand{\nl}{\let\nl\oldnl}}
\makeatother

\newcommand{\mycomment}[1]{\hfill {\footnotesize $\blacktriangleright$ \texttt{#1}}}
\newcommand{\mycommentbeg}[1]{\hfill {\footnotesize $\blacktriangleright$ \texttt{#1}}}
\newcommand{\mycommentend}[1]{\hfill {\footnotesize \texttt{#1}}}

\usepackage{textcomp}
\usepackage{xcolor}
\usepackage{mathrsfs}
\usepackage{bbm}
\usepackage{tcolorbox}

\usepackage{caption}
\usepackage{subcaption}

\makeatletter
\let\NAT@parse\undefined
\makeatother
\usepackage[colorlinks = true, linkcolor = blue, citecolor = blue]{hyperref}

\newcommand{\ldef}{:=}
\newcommand{\rdef}{=:}
\newcommand{\Mcal}[1]{\mathcal{#1}}
\newcommand{\tth}{^{\text{th}}}

\newcommand{\relu}[1]{\left[ #1 \right]_+}

\graphicspath{{./figs/}}

\def \real{\mathbb{R}}
\def \integer{\mathbb{Z}}
\def \one{\mathbf{1}}
\def \x{\mathbf{x}}
\def \neigh{\mathcal{N}}
\def \bneigh{\overline{\mathcal{N}}}

\def \outneigh{\, ^{\text{out}}\neigh}

\def \lag{L}

\def \prob{\textbf{P}}

\def \graph{\mathcal{G}}

\def \flow{\mathcal{F}}
\def \redflow{\hat{\mathcal{F}}}
\def \tempflow{\flow_{\,\text{temp}}}
\def \sflow{\mathcal{H}}
\def \node{\mathcal{V}}
\def \snode{\mathcal{W}}
\def \rednode{\hat{\mathcal{V}}}
\def \tempnode{\node_{\,\text{temp}}}
\def \edg{\mathcal{E}} 
\def \arc{\mathcal{A}}
\def \barc{\overline{\mathcal{A}}}
\def \sarc{\mathcal{Q}}
\def \redarc{\hat{\mathcal{A}}}
\def \temparc{\arc_{\text{temp}}}
\def \simplex{\mathcal{S}}

\def \ne{\mathcal{NE}}

\def \inmat{\mathbf{A}}

\def \pe{p}
\def \lev{u}
\def \dyn{\mathbf{F}}

\def \supp{\mathrm{supp}}

\def \gpath{\Mcal{P}}

\def \M{\mathcal{M}}
\def \z{\mathbf{z}}
\def \a{\mathbf{a}}
\def \d{\mathbf{d}}

\def \rvec{\mathbf{r}}
\def \w{\mathbf{w}}

\def \hill{\Mcal{H}}
\def \bhill{\overline{\Mcal{H}}}
\def \optu{f}
\def \rhom{\rho^\text{max}}
\def \inreach{\, ^{\text{in}}\mathcal{R}}
\def \outreach{\, ^{\text{out}}\mathcal{R}}

\def \bx{\overline{\x}}
\def \brho{\overline{\rho}}
\def \bussd{\overline{U}_{\text{SSD}}(\x^0)}
\def \bunbrd{\overline{U}_{\text{NBRD}}(\x^0)}
\def \bunrpm{\overline{U}_{\text{NRPM}}(\x^0)}
\def \umax{U_{\text{max}}}
\def \umin{U_{\text{min}}}

\def \redflowx{\redflow}
\def \rednodex{\rednode}
\def \redarcx{\redarc}

\def \dumnode{{\node^\prime}}
\def \dumarc{{\arc^\prime}}
\def \dumgraph{{\graph^\prime}}

\def \poslim{\Mcal{L}^+}

\newtheorem{theorem}{Theorem}[section]

\newtheorem{definition}[theorem]{Definition}
\newtheorem{lemma}[theorem]{Lemma}

\newtheorem{remark}[theorem]{Remark}

\newcommand{\thmtitle}[1]{\mbox{}\textit{(#1).}}

\newcommand{\bulletsym}{\hbox{$\bullet$}}
\newcommand{\bulletend}{\relax\ifmmode\else\unskip\hfill\fi\bulletsym}
\newcommand{\squaresym}{\hbox{$\blacksquare$}}
\newcommand{\proofend}{\relax\ifmmode\else\unskip\hfill\fi\squaresym}
\newcommand{\trianglesym}{\hbox{$\blacktriangle$}}
\newcommand{\egend}{\relax\ifmmode\else\unskip\hfill\fi\trianglesym}

\newcommand{\bld}[1]{\mathbf{#1}}

\newenvironment{proofa}
{}{\proofend}

\def \papername{Dynamics of a Stratified Population of Optimum Seeking
  Agents on a Network - Part II : Steady State Analysis}
  \def \shortpapername{Stratified population dynamics on a network: steady state analysis}

\begin{document}


\title{\papername} \author{Nirabhra Mandal and Pavankumar Tallapragada
  \IEEEmembership{Member, IEEE}
  \thanks{
    This work was partially supported by Robert Bosch Centre for
    Cyber-Physical Systems, Indian Institute of Science,
    Bengaluru. N. Mandal was supported by a fellowship grant from the
    Centre for Networked Intelligence (a Cisco CSR initiative) of the
    Indian Institute of Science, Bengaluru. } \thanks{N. Mandal is
    with the Department of Electrical Engineering, Indian Institute of
    Science, Bengaluru. P. Tallapragada is with the Department of
    Electrical Engineering and the Robert Bosch Centre for Cyber
    Physical Systems, Indian Institute of Science, Bengaluru \tt\small
    \{nirabhram, pavant\}@iisc.ac.in } }

\maketitle

\begin{abstract}
  In this second part of our work, we study the steady state of the
  population and the social utility for the three dynamics SSD, NBRD
  and NRPM; which were introduced in the first part. We provide
  sufficient conditions on the network based on a maximum payoff
  density parameter of each node under which there exists a unique
  Nash equilibrium. We then utilize positive correlation properties of
  the dynamics to reduce the flow graph in order to provide an upper
  bound on the steady state social utility. Finally we extend the idea
  behind the sufficient condition for the existence of a unique Nash
  equilibrium to partition the graph appropriately in order to provide
  a lower bound on the steady state social utility. We also illustrate
  interesting cases as well as our results using simulations.
\end{abstract}


\section{Introduction} \label{sec:intro}

\IEEEPARstart{A}{} primary goal in the analysis of evolutionary
dynamics is to characterize the set of equilibrium points and to study
their stability. Many applications might also require the knowledge of
the steady state, given the initial condition of the population, for
higher level control and planning problems. Similarly, knowledge of
the steady state value of the social utility as a function of the
initial population state can be useful for many purposes. In this
paper, we seek to obtain efficiently computable bounds on the steady
state social utility under \emph{stratified smith dynamics} (SSD),
\emph{nodal best response dynamics} (NBRD) and \emph{network
  restricted payoff maximization} (NRPM) as a function of the initial
population state.

\subsection{Literature Survey}

Population games and evolutionary dynamics~\cite{WHS:2010:pop_games}
find application in problems related to distributed control and
formation control \cite{NQ-COM-JBG-GO-AP-EMN:2017:article_app,
  JBG-GO-NQ:2016:similar,
  JMP-GDG-NQ-LPG:2020:ifac_formation_control_app} as well as in social
or socio-technical systems such as transportation and opinion
dynamics. In the context of network games, \cite{EL-CH-MAN:2005:egt1,
  KP-MB-JR-LJS:2015:egt2, LMH-NC-KH:2011, BA-MAN:2014:games_graph,
  BA-GL-YTC-BF-NM-STY-MAN:2017} model a finite population of agents as
nodes and models the interactions using the graph. Other works of
literature \cite{JBG-GO-NQ:2016:similar, GC-FF-LZ:2020:ifac_community,
  JBG-HT:2018, JBG-HT:2018:NEG, JBG-GO-AP-HT:2020:ne_neg_ifac}
consider the nodes of the graph to represent choices with the state of
the population being composed of the fractions of population choosing
different nodes. In these works, the network plays a major role in the
evolution of the population.

In part one of this work\cite{NM-PT:2020:journ1}, we modeled the game
with stratified populations with their choices as the nodes of a
graph, modeled the three dynamics: SSD, NBRD and NRPM. We also
existence and uniqueness of solutions of these dynamics and showed
that their solutions converge to the set of Nash equilibria.

%

\subsection{Contribution}

In this second part of our work, we first provide sufficient
conditions on the graph under which the stratified population game has
a unique Nash equilibrium. Then, for general graphs, we provide a
computationally efficient method for computing bounds on the steady
state social utility. This method is particularly useful in the case
of NBRD and NRPM, as these dynamics rely on an underlying optimization
problem, which makes simulating the full dynamics in order to
determine these steady state values computationally expensive. Indeed,
the method we propose depends only on the initial population state,
the network and the node parameters.

References~\cite{JBG-GO-NQ:2016:similar, JBG-HT:2018, JBG-HT:2018:NEG,
  JBG-GO-AP-HT:2020:ne_neg_ifac} have a setup for the underlying game
closest to ours. However, all these papers assume that Nash
equilibrium is unique and in the relative interior of the
$n$-dimensional probability simplex and results are local. Moreover,
in all these works, the population is not stratified. Our preliminary
work on this topic\cite{NM-PT:2020:ifac} considers only quadratic
cumulative payoff functions as opposed to strictly concave functions
and is concerned only with the convergence analysis. Here, we provide
a sufficient condition on the graph and node parameters under which
the Nash equilibrium is unique and for a general graph provide a
computationally efficient method for computing bounds on the steady
state social utility, given the initial population state.

\subsection{Organization}

The rest of the paper is organized as follows. In Section
\ref{sec:prelim}, we list a few important results from the first part
of this work and outline the problem we wish to address in this
paper. In Section \ref{sec:unique_nash}, we give sufficient conditions
on the graph under which the population game has a unique Nash
equilibrium. In Section \ref{sec:bounds-social-utility} we give
algorithms to provide meaningful bounds on the steady state social
utility. In Section \ref{sec:sims} we provide simulation results and
some interesting examples.

\subsection{Notation and Definitions}

We denote the set of real numbers, the set of non-negative real
numbers and the set of integers with $\real$, $\real_+$ and
$\integer$, respectively. We let $[p,q]_\integer$ be the set of all
integers between $p$ and $q$ (inclusive), \emph{i.e.},
$[p,q]_\integer \ldef \{x \in \integer \,|\, p \leq x \leq
q\}$. $\real^n$ (similarly $\real_+^n$) is the cartesian product of
$\real$ (equivalently $\real_+$) with itself $n$ times. If $\bld{v}$
is a vector in $\real^n$, we denote $\bld{v}_i$ as the $i^\mathrm{th}$
component of $\bld{v}$ and for a vector $v \in \real^n$, we let
$\supp(\bld{v}) \ldef \{i \in [1,n]_\integer \,|\, \bld{v}_i \neq
0\}$. 
We let $\one$ be the vector, of appropriate size, with all its
elements equal to 1 and we let $\bld{e}_i$ be the vector, again of
appropriate size, with its $i^\mathrm{th}$ element equal to 1 and 0
for all other elements. The empty set is denoted by $\varnothing$. If
$\mathcal{Q}$ is an ordered countable set, then $\mathcal{Q}_i$
denotes the $i^\mathrm{th}$ member of $\mathcal{Q}$ and
$|\mathcal{Q}|$ is used to represent the cardinality of
$\mathcal{Q}$. For two sets
$\mathcal{U}, \mathcal{V} \subset \mathcal{Q}$, the set subtraction
operation is denoted by
$\mathcal{U} \setminus \mathcal{V} = \mathcal{U} \cap \mathcal{V}^c$,
where $\mathcal{V}^c$ is the set complement of $\mathcal{V}$ in
$\mathcal{Q}$. $\{i,j\}$ is used to denote an unordered pair while
$(i,j)$ is used to denote an ordered pair. For a vector
$\bld{v} \in \real^n$, $\bld{v} \geq 0$ is used to denote term wise
inequalities. 
By $\relu{.} : \real \to \real_+$ we denote the function that is
defined as $\relu{x} \ldef \max \{x,0\}$. For a function
$f(\x) : \real^n \rightarrow \real$, $\nabla_\x f$ is used to denote
the gradient of $f(.)$ with respect to $\x$, \emph{i.e.}, the $j\tth$
component of $\nabla_\x f$ is $\frac{\partial f}{\partial \x_j}$. We
denote by $\simplex^n_\gamma$, the $n$-dimensional simplex
$\simplex^n_\gamma \ldef \{\bld{v} \in \mathbb{R}^n \,|\,\, \bld{v}
\geq 0 \,\,,\, \one^T \bld{v} = \gamma \}$.

\section{Preliminaries} \label{sec:prelim}%

In this section, we recollect the framework and the essential aspects
about the three dynamics SSD, NBRD, NRPM that we developed
in\cite{NM-PT:2020:journ1}. We let the \emph{population} is composed
of a continuum of \emph{agents} and the network of \emph{choices} is
given by an undirected connected graph $\graph \ldef (\node,\edg)$,
where $\node$ is a set of nodes and
$\edg \subseteq \node \times \node$ is a set of edges. The fraction of
population in node $i$ is represented by $\x_i \in [0, 1]$ and
$\sum_{i \in \node} \x_i = 1$. Thus $\x \in \simplex^{|\node|}_1$
represents the population state.

Let $\pe_i(.) : [0,1] \to \real$ be the function that models the
\emph{cumulative payoff} of the fraction $\x_i$. We assume that the
fraction in each node is stratified and the agents in different strata
of a given node receive different payoffs. Let
$[a,b] \subseteq [0,\x_i]$ be an arbitrary interval. Then the agents
of node $i$ that are in the \emph{strata} $[a,b]$ get an \emph{average
  payoff} of
\begin{equation*}
  \frac{\pe_i(b)-\pe_i(a)}{b-a} .
\end{equation*} 
For a node $i$ if $a \in [0,\x_i]$, then
\begin{equation*}
  \lev_i(a) \ldef \frac{\mathrm{d}\,\pe_i}{\mathrm{d}\,y}(a) \,,
\end{equation*}
is the average payoff that the agents in the strata $[a]$ of node $i$
receive. By \emph{strata} $[a]$ we mean the collection of
\emph{infinitesimal strata} around $a$ in $[0,\x_i]$. We call
$\lev_i(.)$ as the \emph{payoff density function} of node $i$.  We let
$\lev_i(0)$ be the right derivative of $\pe_i(.)$ at \emph{zero} and
$\lev_i(1)$ be the left derivative of $\pe_i(.)$ at \emph{one}. We let
$u(.)$ be the vector whose $i\tth$ element is $u_i(.)$. Through out
this paper, we make the following assumption.
\begin{itemize}
\item[\textbf{(A1)}] For all $i \in \node$, $\pe_i(.)$ is twice
  continuously differentiable and strictly concave. Hence,
  $\forall i \in \node$, $\lev_i(.)$ is a strictly decreasing
  function.
\end{itemize}

The function,
\begin{equation}
	U(\x) \ldef \sum_{i \in \node} \ [\pe_i(\x_i) - \pe_i(0)] \,,
	\label{eq:social_utility}
\end{equation}
which we call as the \emph{social utility function} represents the aggregate payoff of the population as a whole. Note that $U(.)$ is a strictly concave function and $\bld{e}^T_i \nabla U (\x) = u_i(\x_i)$, $\forall \x \in \simplex^N_1$. We let $\neigh^i$ be the set of all neighbors of node $i$ in the graph
$\graph$ and let $\bneigh^{\,i}= \neigh^i \cup \{i\}$. Given the
undirected graph $\graph$, we let
\begin{equation}
  \arc \ldef \bigcup_{\{i,j\} \in \edg} \{ (i,j),(j,i) \} .
  \label{eq:arc_from_edg}
\end{equation}
When an infinitesimal agent in node $i \in \node$ decides to switch to
a node $j \in \neigh^i$, it gets to enter the strata $[\x_j]$ of node
$j$. The dynamics that govern the evolution of the stratified population can be described by the general equation
\begin{equation}
  \dot{\x} = \inmat \,\Delta(\x) \rdef \dyn(\x) \, .
  \label{eq:gen_dyn}
\end{equation} 
Here $\inmat \in \real^{|\node| \times |\arc|}$ is the incidence
matrix of the directed graph $\flow \ldef (\node,\arc)$ and $\Delta(\x) \in \real^{|\arc|}$ is the vector that accumulates the outflows $\delta_{ij}$ along each arc $(i,j) \in \arc$ into a vector. Such a class of dynamics is referred to as \emph{flow balanced dynamics}. The following discussion holds for 
a general case where the total population is $\rho \in \real_+$, \emph{i.e} $\sum_{i \in \node} \x_i = \rho$. 

The sub class of flow balanced dynamics that additionally satisfies the following definition of \emph{strong positive correlation} are referred to as \emph{strongly positively correlated flow balanced dynamics}.
\begin{definition}\thmtitle{Strong positive correlation}\label{def:PC}
  We say that $\dyn(.) : \real^{|\node|} \to \real^{|\node|}$ in \eqref{eq:gen_dyn} is
  strongly positively correlated with $u(.)$ if
  $\,\forall \x \in \simplex^{|\node|}_\rho$, $u_i(\x_i) \geq u_j(\x_j)$ implies
  $\delta_{ij} = 0$, $\forall (i,j) \in \arc$.  \bulletend
\end{definition}
The following lemma enumerates the key result regarding this types of dynamics.

\begin{lemma} \label{lem:ne_fbd}
\thmtitle{Non-emptiness of set of Nash equilibrium for strongly positively correlated dynamics} 
Suppose $\dyn(.)$ is strongly positively
  correlated with $\lev(.)$. Then the set
  \begin{equation}
    \begin{split}
      &\ne^{|\node|}_\rho := \big\{\x \in \simplex^{|\node|}_\rho \,\,\big |\,\, \lev_i(\x_i)
      \geq \lev_j(\x_j),\\
      & \qquad \qquad \qquad \qquad \qquad \forall \, j \in \neigh^i,
      \forall \, i \in \supp(\x)\big\}\,.
    \end{split}
    \label{eq:Nash_eq}
  \end{equation}
  is non-empty. Furthermore, $\ne^{|\node|}_\rho$ is a subset of
  the equilibrium points of $\dyn(.)$ in $\simplex^{|\node|}_\rho$ and
  the unique optimizer of
\begin{equation}
  \begin{split}
    \prob_1(\dumnode,\rho): \quad \optu_\dumnode(\rho) \ldef \ \max_{\x}
    U(\x)\, \text{ s.t. } \x \in \simplex^{|\dumnode|}_\rho \,,
  \end{split}
  \label{eq:P_no_restrict}
\end{equation}
with $\dumnode = \node$ belongs to the set $\ne^{|\node|}_\rho$.  \bulletend
\end{lemma}
The set $\ne^{|\node|}_\rho$ in \eqref{eq:Nash_eq} is the set of \emph{Nash equilibria} of the population game.

The three dynamics, namely, SSD, NBRD and NRPM, which we have
introduced in\cite{NM-PT:2020:journ1} vary in the levels of
coordination displayed by the \emph{infinitesimal agents}. The
$\Delta$ vector is hence calculated in different ways for each of
them. We discuss these dynamics briefly and state key results related
to them.

\begin{itemize}

\item[(a)] \emph{Stratified Smith Dynamics (SSD): }
Here, each agent is selfish and
revises its choice at independent and random time instants by comparing its payoff with the payoff of a neighboring agent. Thus the individual components of $\Delta$ are given by
\begin{equation}
  	\delta_{ij} = \relu{ \lev_j(\x_j) ( \x_i - y_{ij} ) - ( \pe_i(\x_i) -
    \pe_i(y_{ij}) )},
  \label{eq:sd_delta_ij}
\end{equation}
where,
\begin{equation}
	y_{ij} \ldef 
	\begin{cases}
		0, & \text{if } \lev_i^{-1}(\lev_j(\x_j)) < 0 \\
		\lev_i^{-1}(\lev_j(\x_j)), & \text{if } 0 \leq \lev_i^{-1}(\lev_j(\x_j)) \leq \x_i \\
		\x_i, & \text{if } \lev_i^{-1}(\lev_j(\x_j)) > \x_i \, .
	\end{cases}
	\label{eq:yij}
\end{equation}

The key results regarding SSD are are presented in the following theorem.
\end{itemize}

\begin{theorem}\thmtitle{Strong positive correlation of SSD and asymptotic convergence to the set of Nash equilibria} \label{th:ssd_conv}
Suppose that for all $i \in \dumnode$, $\lev_i(.)$ is strictly decreasing. Let the evolution of $\x$ be governed by SSD, dynamics in
  \eqref{eq:gen_dyn} with $\delta_{ij}$ defined in
  \eqref{eq:sd_delta_ij}, with an intial condition
  $\x(0) \in \simplex^{|\node|}_\rho$. Then SSD is strongly
  positively correlated with $u(.)$ and the set of equilibrium
  points of SSD in $\simplex^{|\node|}_\rho$ is the set $\ne^{|\node|}_\rho$ in \eqref{eq:Nash_eq}. Furthermore, as $t \to \infty$, $\x(t)$ approaches $\ne^{|\node|}_\rho$
  defined in \eqref{eq:Nash_eq} and $U(\x(t))$ converges to a constant. \bulletend
\end{theorem}

\begin{itemize}

\item[(b)] \emph{Nodal Best Response Dynamics (NBRD): }
Here, the agents coordinate with each other at the nodal level. The fraction of population in each node decides to redistribute itself among its neighbors by computing its best
response to the current configuration of $\x$; while assuming that the fraction in the neighboring nodes do not change. Thus $\delta_{ij}$ for NBRD is computed from the following optimization problem.
\begin{align}
  \label{eq:p_nbrd-prob}
  \prob^{\,i}_2:\\
  \notag \max_{\{\d_{ij} | j \in \bneigh^{\,i}\}}%
  & \sum_{j \in \neigh^i}[\pe_j(\x_j+\d_{ij}) - \pe_j(\x_j)] \\
  &\notag \qquad \qquad \qquad \quad  + [\pe_i (\d_{ii}) - \pe_i(0)]\\
  \notag \text{s.t. } & \d_{ii} + \sum_{j \in \neigh^i}
                        \d_{ij} = \x_i, \quad \d_{ij} \geq 0, 
                        \ \forall j \in \neigh^{i}. 
\end{align}
Thus if $\{\d^*_{ij}\}_{(i,j) \in \arc}$ is the subset of optimizers of $\{\prob^{\,i}_2\}_{i \in \node}$, then $\delta_{ij} = \d^*_{ij}$ $\forall (i,j) \in \arc$. The key results related to NBRD are listed below.
\end{itemize}

\begin{theorem}\label{th:NBRD_full}\thmtitle{Strong positive correlation; existence and uniqueness of solutions and asymptotic convergence to the set of Nash equilibria for NBRD} 
Let $\{\d^*_{ij}\}_{(i,j) \in \arc}$ be the subset of optimizers of
  $\{\prob^{\,i}_2\}_{i \in \node}$. Let NBRD be the dynamics in \eqref{eq:gen_dyn} with
  $\delta_{ij} = \d^*_{ij}$ $\forall (i,j) \in \arc$. NBRD is strongly positively correlated with $\lev(.)$. For each $\x(0) \in \simplex^{|\node|}_\rho$, NBRD has a unique
  solution that exists for $\forall \, t \geq 0$. The set of
  equilibrium points of NBRD in $\simplex^{|\node|}_\rho$ is $\ne^{|\node|}_\rho$ defined in
  \eqref{eq:Nash_eq}. Further, if $\x(t)$ evolves according to NBRD, then as $t \to \infty$, $U(\x(t))$ converges to a constant and $\x(t)$ approaches $\ne^{|\node|}_\rho$. \bulletend
\end{theorem}

\begin{itemize}
\item[(c)] \emph{Network Restricted Payoff Maximization (NRPM): }
Here, the agents coordinate across the entire population
and evolve under a centralized decision scheme. They do so in order to
maximize the social utility of the entire population.
\begin{subequations}
  \label{eq:prob_nrpm}
  \begin{align}
    &\prob_3: \notag\\
    &\max_{\z, \d} \sum_{i \in
      \node}[\pe_i(\z_i ) - \pe_i(0)] \quad = \quad \max_{\z, \d} U(\z)
    \\
    &\text{s.t. } \z_i = \x_i + \sum_{j \in \bneigh^{\,i}} (\d_{ji} -
      \d_{ij}) , \ \forall \, i \in
      \node, \label{eq:nrpm-flow-balance}\\
    &\sum_{j \in \bneigh^{\, i}}\d_{ij}  = \x_i, \ \forall
      \, i \in \node, \label{eq:nrpm-xi-reallocation}\\
    &\d_{ij} \geq 0, \ \forall \, (i,j) \in
      \barc \ldef \arc \cup \{(i,i) \ | \ i \in \node \}
      . \label{eq:nrpm-non-neg}
  \end{align}
\end{subequations}

Thus, if $(\z^*,\d^*)$ is an optimizer of $\prob_3$, then the dynamics is given by 
\begin{equation}
  \label{eq:NRPM_ODE}
  \dot{\x} = \inmat \,\Delta(\x) = \z^*(\x) - \x .
\end{equation}
Key results regarding NRPM are listed next.
\end{itemize}
\begin{lemma}\label{lem:p3_unique}
\thmtitle{The resultant node fractions in any
    optimal solution of $\prob_3$ are unique}
  Consider the problem $\prob_3$ in \eqref{eq:prob_nrpm} with
  $\x \in \simplex^{|\node|}_\rho$. Let $\mathcal{OP}$ be the set of all optimizers
  of $\prob_3$. Then $\forall \, (\z^*,\d^*) \in \mathcal{OP}$,
  $\bld{z}^* = \x + \inmat\,\d^*$ is unique. \bulletend  
\end{lemma}

\begin{lemma}\thmtitle{Properties of optimizers of
    $\prob_3$} \label{lem:p3_opt}
  Suppose $(\z^*,\d^*)$ is an optimizer of $\prob_3$. Then
\begin{equation}
	\lev_j(\z^*_j) \geq \lev_i(\z^*_i), \,\, \forall (i,j) \in \arc \,\,\text{s.t.}\,\, \d^*_{ij} > 0 \,.
	\label{eq:level_order_nrpm}
\end{equation}
\bulletend
\end{lemma}

\begin{theorem}\thmtitle{Existence and uniqueness of solutions of NRPM and asymptotic convergence to the set of Nash equilibria}
  Let NRPM be the dynamics in~\eqref{eq:NRPM_ODE}. For each $\x(0) \in \simplex^{|\node|}_\rho$, NRPM has a unique
  solution that exists for $\forall \, t \geq 0$. The set of
  equilibrium points of NRPM in $\simplex^{|\node|}_\rho$ is $\ne^{|\node|}_\rho$ defined in
  \eqref{eq:Nash_eq}. Further, if $\x$ evolves according to NRPM, then as $t \to \infty$, $U(\x(t))$ converges to a constant and $\x(t)$ approaches $\ne^{|\node|}_\rho$. 
  
\bulletend
  \label{th:PM_conv}
\end{theorem}

We provide another definition and a remark next to conclude the
preliminaries.
\begin{definition}\label{def:nash_conv_dyn}
\thmtitle{Nash convergent dynamics}
Let $\dumgraph = (\dumnode,\dumarc)$ be the underlying graph.
Any flow balanced dynamics $\dyn(.)$ on $\dumgraph$ that asymptotically approaches set of Nash equilibria $\ne^{|\dumnode|}_\rho$ from any initial condition in $\simplex^{|\dumnode|}_\rho$ is termed as Nash convergent dynamics. \bulletend
\end{definition}

\begin{remark}\label{rem:nash_conv}
\thmtitle{SSD, NBRD and NRPM are Nash Convergent}
By Theorems \ref{th:ssd_conv}, \ref{th:NBRD_full} and \ref{th:PM_conv}; SSD, NBRD and NRPM are Nash convergent on $\graph$.
\bulletend
\end{remark}

\subsection{Problem Setup} \label{sec:setup}%

In general SSD, NBRD and NRPM may not converge to the same state. In
this paper, we seek to study the steady state and the steady state
social utility with the knowledge of the structure of the graph
$\graph = (\node,\edg)$ and the initial condition, which we denote by
$\x^0 \ldef \x(0)$.

From Remark \ref{rem:nash_conv} SSD, NBRD and NRPM are all Nash convergent dynamics.
Furthermore, for all three dynamics, the social utility converges to a constant value. We define $\bussd$, $\bunbrd$ and $\bunrpm$ to be the
values to which the social utility converges to under the three different
dynamics from an initial condition $\x^0$.

In Section~\ref{sec:unique_nash}, we provide a sufficient
condition under which the steady state state of any Nash convergent dynamics is the same. Then in Section~\ref{sec:bounds-social-utility}, we
provide bounds on the steady state social utility for the three
dynamics.
 
\section{Unique Nash Equilibrium}\label{sec:unique_nash}

In this section, we give sufficient conditions on the graph $\graph$
and the payoff density functions, under which the set of Nash
equilibria is a singleton. We define for each node $i \in \node$, the
\emph{maximum payoff density parameter} (MPDP) as
\begin{equation}
  \a_i \ldef \max_{y \in [0,1]} \lev_i(y) = \lev_i(0) ,
  \label{eq:a}
\end{equation}
where we use the fact that $\lev_i(.)$ is a strictly decreasing
function in $[0,1]$ for all $i \in \node$.

We formally define a path in a graph next.

\begin{definition}\thmtitle{Path in graph $\graph = (\node,\edg)$}
  A path, in graph $\graph = (\node,\edg)$, between $i \in \node$ and
  $j \in \node$ is given by an ordered set of non-repeating elements
  $\gpath(i,j) \subseteq \node$ such that $\gpath_1(i,j) = i$,
  $\gpath_{n}(i,j) = j$ and
  $\{\gpath_k(i,j),\gpath_{k+1}(i,j)\} \in \edg$,
  $\forall k \in {[1,n-1]_\integer}$, where
  $\gpath_k(i,j)$ denotes the $k\tth$ element in
  $\gpath(i,j)$ and $n = |\gpath(i,j)|$.
  \bulletend
\end{definition}
Using this definition and \eqref{eq:a}, we define a \emph{path with quasi-concave MPDP's} next.
\begin{definition}\label{def:path_qc}
  \thmtitle{Path with quasi-concave MPDP's} %
  Suppose $\gpath(i,j)$ is a path between $i \in \node$ and
  $j \in \node$ and let $n \ldef |\gpath(i,j)|$. Let
  $\pi(k) \ldef \gpath_k(i,j)$. We say $\gpath(i,j)$ is a path with
  quasi-concave MPDP's if and only if
  $\forall \, k,l \in [1,n]_\integer$ such that $l \geq k$ and
  $\forall m \in [k,l]_\integer$,
  $\displaystyle \a_{\pi(m)} \geq \min \{\a_{\pi(k)},\a_{\pi(l)}\}$.

\bulletend
\end{definition}
Definition \ref{def:path_qc} can be interpreted in the following way. For a path with quasi-concave MPDP's, if we arrange the MPDP's in order of the nodes visited in the path, then the MPDP's form a quasi-concave function. Using this interpretation and Definition \ref{def:path_qc}, the following result can be immediately stated regarding the nature of MPDP's in such a path.

\begin{lemma}\thmtitle{Monotonicity of MPDP's in a path with quasi-concave MPDP's} \label{lem:monotone_path}
  Let
  $\gpath(i,j)$ be a path with quasi-concave MPDP's between $i$ and
  $j$, where $i,j \in \node$. Let $\pi(k) \ldef
  \gpath_k(i,j)$ and $n =
  |\gpath(i,j)|$. Then $\forall \, k \in [1,n]_\integer$, either
  $\a_i  = \a_{\pi(1)} \leq \a_{\pi(2)} \leq \cdots \a_{\pi(k)}$ or
  $\a_{\pi(k)} \geq \cdots \geq \a_{\pi(n-1)} \geq \a_{\pi(n)} = \a_j$.
  \bulletend
\end{lemma}
A proof of this is provided in Appendix \ref{proof:monotone_path}.

Using Definition \ref{def:path_qc}, we define a \emph{quasi-concave hill} next.

\begin{definition}\label{def:qch}
\thmtitle{Quasi-concave hill or QCH} 
We call $\graph = (\node,\edg)$ to be a quasi-concave hill or QCH if and only if there exists a path with quasi-concave MPDP's between every two distinct nodes (\emph{i.e.} $i,j \in \node$ and $i \neq j$).  \bulletend
\end{definition}
Note that Definition \ref{def:qch} requires the graph $\graph$ to be
connected. The `hill' in the name of a QCH is given to evoke the
interpretation that the agents are always seeking to climb a hill to
increase their payoffs. The following lemma shows that in a QCH, if
the population is in an equilibrium (in the set of Nash equilibria),
then the payoff densities across the non-empty nodes is the same.

\begin{lemma} \label{lem:same_lev} \thmtitle{Same payoff density
    across non-empty nodes in a population at equilibrium in a QCH} %
  Let $\graph$ be a QCH and let $\x \in \ne^{|\node|}_\rho$. Then
  $\forall \, i,j \in \supp(\x)$, $\lev_i(\x_i) = \lev_j(\x_j)$.
  
  \bulletend
\end{lemma}
We present the proof in Appendix \ref{proof:same_lev}.

Finally, in this section, we show that for a QCH the Nash equilibrium
is unique. We illustrate the consequence of this fact in the remark
following the theorem, proof of which is presented in Appendix \ref{proof:unique_nash}
\begin{theorem} \label{th:unique_nash} \thmtitle{Uniqueness of Nash
    equilibrium for a QCH} %
  Suppose $\rho \geq 0$. If $\graph = (\node,\edg)$ is a QCH then
  $\ne^{|\node|}_\rho$ in \eqref{eq:Nash_eq} is a singleton and the
  unique $\x \in \ne^{|\node|}_\rho$ is the unique optimizer of
  $\prob_1$ in \eqref{eq:P_no_restrict}.
  \bulletend
\end{theorem}

\begin{remark}\label{rem:same_conv_st}
  \thmtitle{Same steady state state of Nash convergent dynamics for QCH}
  The immediate consequence of Theorem \ref{th:unique_nash} is that
  $\forall \x(0) \in \simplex^{|\node|}_\rho$ the steady state of any
  Nash convergent dynamics is the same. Thus, from
  Remark~\ref{rem:nash_conv}, we can say that if the population
  dynamics is governed by any one of SSD, NBRD and NRPM then
  $\forall \x(0) \in \simplex^{|\node|}_\rho$ the population converges
  to the unique Nash equilibrium.  \bulletend
\end{remark}

The utility of Theorem~\ref{th:unique_nash} and
Remark~\ref{rem:same_conv_st} extends to the case in which $\graph$ is
not a QCH. In particular, we can use them to provide bounds on the
steady state social utility of SSD, NBRD and NRPM. This is the topic
of the next section.

\section{Bounds on the Steady State Social Utility}
\label{sec:bounds-social-utility}

In general, the steady state for SSD, NBRD and NRPM may not be the
same. The graph $\graph$ being a QCH is a sufficient condition for all
these dynamics to share a common steady state irrespective of the
initial condition in $\simplex^{|\node|}_\rho$. In this section, we
provide meaningful bounds on the steady state social utility of SSD,
NBRD and NRPM in the case where $\graph$ is not a QCH. In the
construction of these bounds, the results stated in
Section~\ref{sec:unique_nash} serve as important building blocks.

Recall that $\graph = (\node,\edg)$ is the underlying undirected graph
and $\arc$ the associated arc set where for each undirected edge
$\{i, j\} \in \edg$, we have directed arcs $(i,j)$ and $(j,i)$ in
$\arc$. Moreover, $\arc$ contains no other arcs. Based on the payoff
functions $\lev_i(.)$, it is possible to reduce $\arc$ based on the
fact that some $\delta_{ij}$'s will remain \emph{zero} for SSD, NBRD
and NRPM no matter what the initial condition is. We characterize this
in the following lemma.
\begin{lemma}\label{lem:arc_remove}
\thmtitle{Reduction of flow graph for SSD, NBRD and NRPM}
Let $\lev_i(.)$'s be strictly decreasing functions. Let $\theta \in \real^{|\node|}$ be a vector with individual components
  $\theta_i \geq 0$. If $\x_i \in [0, \theta_i]$, then
\begin{align}
	\delta_{ij}(\x) = 0, \ \forall (i,j) \in \arc \text{ s.t. } \lev_i(\theta_i) \geq \lev_j(0) \,,
	\label{eq:delta_zero_cond}
\end{align}
for SSD and NBRD. Further if $\z^*_i(\x) \in [0, \theta_i]$,
then~\eqref{eq:delta_zero_cond} holds for NRPM as well.
\bulletend
\end{lemma}
A proof for the same is provided in Appendix \ref{proof:arc_remove}.
Thus, for a given initial state $\x(0) = \x^0$, if we are able to find
a uniform (in time) upper bound $\theta_i$ on $\x_i$ and $\z^*_i(\x)$
for each $i \in \node$ then for the arcs $(i,j) \in \arc$ with the
property $\lev_i(\theta_i) \geq \lev_j(0)$, the flow $\delta_{ij} = 0$
for all time $t \geq 0$ and hence the arc $(i,j)$ can be removed
without modifying the evolution of $\x$. Now as $\forall i \in \node$,
$\x_i \in [0,1]$; $1 \geq \theta_i$, $\forall i \in \node$. Thus
\emph{one} is a good initial estimate of $\theta_i$. Using this idea we can remove the arcs in the set
\begin{equation}
	\{(i,j)\in\arc \,|\, \lev_i(1) \geq \lev_j(0)\}
	\label{eq:red_arc}
\end{equation} 
without affecting the dynamics. We call this
reduction as a \emph{static reduction} of the arc set as this
reduction is carried out solely based on the properties of $\graph$
and the functions $\pe_i(.)$. Note that atmost one arc per pair of
nodes is removed in this process as $\lev_i(1) \leq \lev_j(0)$ and
$\lev_j(1) \leq \lev_i(0)$ are not simultaneously possible due to the strictly decreasing nature of $\lev_i(.)$'s.

In order to reduce the graph further, the initial condition $\x^0$ must be taken into account. In fact, only the nodes which have a directed path in $\redflow$ from some node in $\supp(\x^0)$ are the ones that may participate in the evolution of $\x$. The following remark describes how Algorithm \ref{algo:reduce_graph} utilizes Lemma \ref{lem:arc_remove} to carry out such a reduction.

\begin{algorithm}
\caption{\texttt{reduceGraph}$(\flow,\x^0)$}
\label{algo:reduce_graph}
\SetAlgoLined
\KwData{$\flow = (\node,\arc), \x^0$}
\KwResult{$\redflowx$}


$\redflowx \gets (\node,\arc)$\\

$\tempflow \gets (\varnothing,\varnothing)$\\

\While{ \emph{$\redflowx \neq \tempflow$}}{
	$\tempflow \gets \redflowx$\\
	
	$\redflowx \gets$ \texttt{repeatReduction} $(\redflowx,\x^0)$
}

\textbf{return:} $\redflowx$
\end{algorithm}

\begin{algorithm}
\caption{\texttt{repeatReduction}$(\flow,\x^0)$}
\label{algo:repeat_reduction}
\SetAlgoLined
\KwData{$\flow = (\node,\arc), \x^0$}
\KwResult{Reduced graph}


\For{\emph{$i \in \node$}}{ %
  $\Mcal{R}^i \gets \{ j \in \supp(\x^0) \,|\, \exists$ a directed
  path \label{stp:theta_beg}\\
  \nonl $\qquad \qquad \qquad \qquad \quad \,\,$ in $\flow$ from $j$ to $i\}$\\
  $\theta_i \gets \sum_{j \in \Mcal{R}^i} \x_j$ \label{stp:theta_end}\\
  \nonl \mycommentbeg{estimates the maximum population}\\
  \nonl \mycommentend{fraction that can visit $i$}\\
}
$\temparc \gets \arc \setminus \{(i,j) \in \arc \,|\, \lev_i(\theta_i) \geq \lev_j(0) \}$\\
$\tempflow \gets (\node,\temparc)$\\
$\tempnode \gets \supp(\x^0)$ \label{stp:start_red} \\
\nonl\mycommentbeg{$\supp(\x^0)$ is always a subset of}\\
\nonl\mycommentend{the reduced node set}\\

$\tempnode \gets \tempnode \cup \{i \in \node \,|\, \exists \text{ a directed path in } \tempflow$\\
\nonl $\qquad \qquad \qquad \qquad  \quad \,\, \text{ from some } j \in \supp(\x^0) \text{ to } i\}$\\

$\temparc \gets \temparc \setminus \{(i,j) \in \temparc \,|\, i \notin \tempnode$\\
\nonl $ \qquad \qquad \qquad \qquad \qquad \qquad \quad \text{ or } j \notin \tempnode\}$\label{stp:end_red}\\
\nonl \mycomment{remove hanging arcs}\\

\textbf{return:} $\tempflow = (\tempnode,\temparc)$
\end{algorithm}

\begin{remark}\label{rem:graph_red_init_st}
  \thmtitle{Graph reduction using initial state} The function
  $\texttt{reduceGraph}()$ described in Algorithm
  \ref{algo:reduce_graph} takes in the graph $\flow$ and $\x^0$ as
  inputs and returns the reduced graph
  $\redflowx = (\rednodex, \redarcx)$. Algorithm
  \ref{algo:repeat_reduction} simultaneously estimates $\theta_i$ and
  removes nodes and arcs of $\flow$ until the graph cannot be reduced
  further. Steps \ref{stp:theta_beg}-\ref{stp:theta_end} are used to
  estimate $\theta_i$ by setting it as the sum of all possible
  population fractions that can reach $i$, respecting the structure of
  the graph and the node parameters. Then Steps
  \ref{stp:start_red}-\ref{stp:end_red} only consider $\supp(\x^0)$
  and those nodes in $\flow$ that have a directed path from some node
  in $\supp(\x^0)$. The other nodes and hanging arcs are removed. Note
  that the first pass of Algorithm \ref{algo:reduce_graph} performs
  the static reduction by removing the arcs in the set in
  \eqref{eq:red_arc}. This is because $\flow$ is strongly connected
  and hence the $\theta$ estimated in the first pass will be
  $\rho\one$. Algorithm \ref{algo:reduce_graph}, thus, returns the
  reduced graph $\redflowx = (\rednodex, \redarcx)$ with the following
  properties:
  \begin{itemize}
  \item $\supp(\x^0)$ is contained in $\rednodex$;
  \item $\redarcx$ does not contain arcs $(i,j)$ of $\arc$
    with the property $\lev_i(\theta_i) \geq \lev_j(0)$ for the
    estimated $\theta$.
  \end{itemize}
  $\redflowx$ is called the \emph{initial condition reduced graph} or
  ICRG of $\graph$.

Now, the fact that the Steps \ref{stp:theta_beg}-\ref{stp:theta_end} of Algorithm \ref{algo:repeat_reduction} refine, with each pass, the upper bound on both $\x(t)$
  and $\z^*(\x(t))$ can be easily seen from the fact that $\z^*(\x)$
  is just a simple rearrangement of the population state $\x$ among
  $\node$. Thus by Lemma \ref{lem:arc_remove}, the
  evolution of the population under $\graph$ is equivalent to the
  evolution of the population under $\redflowx$.
\bulletend
\end{remark}

The ICRG is useful in providing reasonable upper and lower bounds on the steady state social utility. We demonstrate this in the next two sections.

\subsection{Upper Bound on Steady State Social Utility} \label{sec:upper_bound}

In order to compute an upper bound on the steady state social utility,
we try to compute the best social utility that the population as a
whole can receive starting from $\x^0$. This best social utility is in
the space of all possible dynamics that have the property that the
evolution is same under $\graph$ and $\redflowx$. From
Remark~\ref{rem:graph_red_init_st}, we see that $\theta_i$ provided by
\texttt{reduceGraph}$(\flow,\x^0)$, Algorithm~\ref{algo:reduce_graph},
is an upper bound on $\x_i(t)$ and $\z^*_i(\x(t))$ for all $t \geq
0$. Thus, if $\redflowx$ is the ICRG of $\graph$ and
$(i,j) \notin \redarcx$ then by Lemma~\ref{lem:arc_remove}
$\delta_{ij}(\x(t)) = 0$ for all $t \geq 0$ in the case of SSD, NBRD
and NRPM. Moreover if $i \in \graph$ but $i \notin \redflowx$, then
$\x_i(t) = 0$, $\forall t \geq 0$. For a node $i \in \node$ of
$\graph$, we define the set of in-reachable nodes of $i$ as
$\inreach^i \ldef \{j \in \rednodex \,|\, \exists \text{ a directed
  path from } j \text{ to } i \text{ in } \redflowx\} \cup \{i\}$ and
the set of out-reachable nodes of $i$ as
$\outreach^i \ldef \{j \in \rednodex \,|\, \exists \text{ a directed
  path from } i \text{ to } j \text{ in } \redflowx\} \cup \{i\}$. We
use $\rvec_{ij}$ to represent the outflow from $i$ to
$j \in \outreach^i$ and let $\rvec$ be the vector that accumulates all
the $\rvec_{ij}$ into a vector. Then the optimum value of the
optimization problem $\prob_4$ can be used as an upper bound on the
steady state utility.
\begin{equation}
  \begin{split}
    \prob_4: \qquad &
    \begin{split}
      U_\text{max} \ldef \max_{\w,\rvec} & \,\, \sum^{|\node|}_{i=1} \left[ \pe_i \left( \w_i \right) - \pe_i(0) \right] \\
      \text{s.t.} & \,\, \w_i = \sum_{j \in \inreach^i} \rvec_{ji}, \ \forall i \in \node\\
      &\, \,\,\,\, \rvec_{ij} \geq 0 \ \forall i \in \node, \forall j \in \outreach^i\\
      & \,\, \sum_{j \in \outreach^i} \rvec_{ij} = \x^0_i, \ \forall i
      \in \node \, .
    \end{split}
  \end{split}
  \label{eq:u_max_opt_prob}
\end{equation}
Note that $\prob_4$ is visually similar to $\prob_3$ of NRPM but they
are very different problems. In $\prob_3$, the decision variables in
$\d$ restricts the movement of the population fraction in each node to
itself and among its neighbors in $\graph$, \emph{i.e.} $\d_{ij}$
captures the outflow form node $i \in \node$ to node $j \in \neigh^i$
(the neighbor set only). In $\prob_4$, on the other hand, the decision
variables in $\rvec$ allows the movement of the population fraction
from each node to any node that is out-reachable from it in
$\redflowx$. Thus, while $\prob_3$ gives us the instantaneous flows in
NRPM, $\prob_4$ gives us the socially optimal longterm redistribution
of the population starting from $\x^0$ and under the path constraints
imposed by $\redflowx$.

We present the main result of this section next and show the proof in Appendix \ref{proof:upper_bound}.
\begin{theorem}\label{th:upper_bound}
\thmtitle{Upper bound on steady state social utility}
Suppose $\x^0 \in \simplex^{|\node|}_\rho$. Then
\begin{align*}
	\bussd \leq \umax; \,\, \bunbrd \leq \umax; \,\, \bunrpm \leq \umax \,.
\end{align*}
\bulletend
\end{theorem}

In the next subsection, we provide algorithms to compute lower bounds
on the steady state social utility.

\subsection{Lower Bound on Steady State Social Utility}

In order to provide a meaningful lower bound on the steady state social
utility, we partition the graph into certain subgrpahs, each of which
is a directed QCH, which we define formally in the sequel. Then, using
the properties in Section~\ref{sec:unique_nash}, we utilize
\eqref{eq:P_no_restrict} to compute the social utility of the
population among each group independently for different allocations of
population fractions in the group. Then, the problem is converted to
one of optimal allocations to the collection of directed QCHs so that
social utility is minimized.

We now formally define a directed QCH and subsequently, we also define
other useful definitions. \textbf{Note:} up to the very end of this
subsection, we present in the context of an arbitrary graph and in the
end we apply it to ICRG $\redflowx$.
\begin{definition}\thmtitle{Directed quasi-concave hill}
  We call a di-graph to be a \emph{directed quasi-concave hill} (DQCH)
  if and only if its corresponding undirected graph is a
  QCH. \bulletend
\end{definition}

\begin{definition}\thmtitle{Directed quasi-concave hill component or
    DQCH component}
  A subgraph $\sflow = (\snode, \sarc)$ of the graph
  $\flow = (\node, \arc)$ is said to be a \emph{directed quasi-concave
    hill component} or \emph{DQCH component} if and only if $\sflow$
  is a DQCH. 
  
  \bulletend
\end{definition}

%
Each DQCH component can be classified into one of two types:
\emph{attractive} or \emph{non-attractive}. We define these next and
also describe their significance.

\begin{definition}\thmtitle{Attractive DQCH component or A-DQCH component}
  Let $\flow = (\node, \arc)$ be a directed graph and let
  $\outneigh^i(\flow)$ be the set of out-neighbors of a node
  $i \in \node$ in $\flow$.  A subgraph $\sflow = (\snode, \sarc)$ of
  the graph $\flow$ is said to be an \emph{attractive directed
    quasi-concave hill component} or \emph{A-DQCH component} of
  $\flow$ if and only if
\begin{itemize}
\item $\sflow$ is a DQCH;
\item $\forall i \in \snode$, $\outneigh^i(\flow) \subseteq \snode$
  and $\sarc = \{(i,j) \,|\, j \in \outneigh^i(\flow)\}$.
\end{itemize}
\end{definition}

\begin{definition}\thmtitle{Non attractive DQCH component or NA-DQCH component}
  Let $\flow = (\node, \arc)$ be a directed graph and let
  $\outneigh^i(\flow)$ be the set of out-neighbors of a node
  $i \in \node$ in $\flow$. A subgraph $\sflow = (\snode, \sarc)$ of
  the graph $\flow$ is said to be a \emph{non-attractive directed
    quasi-concave hill component} or \emph{NA-DQCH component} of
  $\flow$ if and only if
  \begin{itemize}
  \item $\sflow$ is a DQCH;
  \item $\exists \, i \in \snode$ and $j \in \outneigh^i(\flow)$ such
    that $j \notin \snode$. \bulletend
  \end{itemize}
\end{definition}
From the definition of an A-DQCH component, it is clear that if some
fraction of population starts inside it, then the fraction remains
there forever as there are no outgoing arcs from such a
component. This is the main reasoning for the nomenclature. Moreover,
for every DQCH component, if the fraction of population in the said
component is known apriori, we can use \eqref{eq:P_no_restrict} to
calculate the steady state social utility for that population
fraction. Thus, the graph $\redflowx$ can be partitioned into such
DQCH components in order to formulate a suitable optimization problem
(like in Section \ref{sec:upper_bound}) for providing a lower bound on
the steady state social utility. Note that such a partition is not
unique. We use the following partition to tighten the lower bound as
much as possible.

\begin{definition}\label{def:mac_partition}
\thmtitle{Maximal attractive component partition or MAC partition}
  A collection of subgraphs
  $\{\hill^q = (\snode^q, \sarc^q)\}_{q \in [1,n]_\integer}$ is said
  to be a \emph{maximal attractive component partition} or \emph{MAC
    partition} of $\flow = (\node,\arc)$ if and only if
\begin{itemize}
\item $\forall q, r \in [1,n]_\integer$ such that $q \neq r$,
  $\snode^q \cap \snode^r = \varnothing$,
  $\sarc^q \cap \sarc^r = \varnothing$;
\item $\bigcup_{q \in [1,n]_\integer} \snode^q = \node$,
  $\bigcup_{q \in [1,n]_\integer} \sarc^q \subseteq \arc$;
\item $\forall q \in [1,n]_\integer$, $\hill^q$ is a DQCH component of
  $\flow$;
\item $\forall q \in [1,n]_\integer$ if $\hill^q$ is an NA-DQCH
  component, then $\hill^q$ does not contain any A-DQCH component of;
\item $\forall q \in [1,n]_\integer$ if $\hill^q$ is an A-DQCH
  component, then there does not exist any
  $\Mcal{I} \subseteq [1,n]_\integer$ such that
  $\bhill^{\,q} \ldef \big( \bigcup_{r \in \Mcal{I} \cup \{q\}}
  \snode^r, \bigcup_{r \in \Mcal{I} \cup \{q\}} \sarc^r \big)$ is an
  A-DQCH component. \bulletend
\end{itemize}
\end{definition}

Note that a MAC partition of a directed graph always exists as a node
in itself is a DQCH component of a graph. So is a combination of two
neighboring nodes of the graph. Thus, a way to find a MAC partition
would be to locate all the A-DQCH components of $\flow$ and then club
the remaining nodes and arcs into different NA-DQCH components. This
is the main logic behind Algorithm \ref{algo:mac_partition} and the
supporting Algorithm \ref{algo:find_qch_comp}. In particular, the
function $\texttt{MACPartition}()$, which we describe in Algorithm
\ref{algo:mac_partition}, takes $\flow$ as an input and returns a MAC
partition of the same.

\begin{algorithm}
\caption{\texttt{MACPartition}$(\flow,\a)$}
\label{algo:mac_partition}
\SetAlgoLined \KwData{$\flow = (\node,\arc), \a$} \KwResult{MAC
  Partition $\hill$}
$\hill \gets$ \texttt{makeQCHComp}$(\flow,\a)$ \label{stp:get_qch_comp}\\
$\node_\text{regroup} \gets \varnothing$\\
$\arc_\text{regroup} \gets \varnothing$\\
\For{\emph{$\hill^q \in \hill$}}{
  \nonl \mycomment{$\hill^q = (\snode^q,\sarc^q)$}\\
  tempVar $\gets \varnothing$\\
  \While{\emph{$\snode^q \neq$ tempVar}}{
    tempVar $\gets \snode^q$\\
    \For{\emph{$i \in \snode^q$}}{
      \If{\emph{$\exists \, j \in \outneigh^i(\flow)$ such that
          $j \notin \snode^q$}}{
        \nonl \mycomment{$\outneigh^i(\flow)$ is with respect to $\flow$}\\
        $\snode^q \gets \snode^q \setminus
        \{i\}$ \label{stp:remove_node_adqch}\\
        \nonl \mycommentbeg{remove nodes with out-neighbors}\\
        \nonl \mycommentend{not in $\snode^q$ and try to}\\
        \nonl \mycommentend{ make $\snode^q$ an A-DQCH}\\
        $\sarc^q \gets \sarc^q \setminus \{(i,k) \in \sarc^q\} \setminus \{(k,i) \in \sarc^q\}$\\
        \nonl \mycomment{remove corresponding arcs}\\
        $\node_\text{regroup} \gets \node_\text{regroup} \cup \{i\}$ \label{stp:v_regroup}\\
        $\arc_\text{regroup} \gets \arc_\text{regroup} \cup \{(i,k) \in \arc\}$ \label{stp:a_regroup}\\
        Remove $\a_i$ from $\a$\\
			}
		}
	}\label{stp:end_for_hq}
}
$\flow_\text{regroup} \gets (\node_\text{regroup},\arc_\text{regroup})$\\
$\hill \gets \hill \,\,\cup$ \texttt{makeQCHComp}$(\flow_\text{regroup},\a)$ \label{stp:regroup}\\
\nonl \mycommentbeg{regroup the remaining nodes and arcs}\\
\nonl \mycommentend{into NA-DQCH components}\\
\textbf{return:} $\hill$
\end{algorithm}

\begin{algorithm}
\caption{\texttt{makeQCHComp}$(\flow,\a)$}
\label{algo:find_qch_comp}
\SetAlgoLined
\KwData{$\flow = (\node,\arc), \a$}
\KwResult{QCH components}
visitSet $\gets \varnothing$\\
\nonl \mycomment{keeps track of nodes already visited}\\
$\hill \gets \varnothing$\\
\While{\emph{visitSet $\neq \node$}}{
	$i \gets \arg\max(\a)$\\
	$\tempnode \gets \{ i \}$ \label{stp:highest_node}\\
	\nonl \mycommentbeg{node with highest MPDP is}\\
	\nonl \mycommentend{definitely in a QCH}\\
	$\tempnode \gets \tempnode \cup \{ j \in \node \,\,|\,\, \exists \text{  a path with}$ \label{stp:mpdp_path}\\
	\nonl $\qquad \text{ quasi-concave MPDP's between } j \text{ and } i \text{ in } \flow\}$\\
	\nonl \mycommentbeg{these consider paths in the}\\
	\nonl \mycommentend{undirected sense}\\
	$\temparc \gets \{(i,j) \in \arc \,\,|\,\, i,j \in \tempnode\}$\\
	visitSet $\gets$ visitSet $\cup \,\,\tempnode$\\
	Remove indices $\{i \in \text{visitSet} \}$ from $\a$ \\
	$\hill \gets \hill \cup \{(\tempnode,\temparc)\}$
}
\textbf{return:} $\hill$
\end{algorithm}

\begin{remark} \label{rem:compute_macp}
\thmtitle{Computation of MAC partitions of $\redflowx$}
The function \texttt{MACPartition}$(\redflowx,\a)$ in Algorithm \ref{algo:mac_partition} computes a MAC partition of $\redflowx$. The function \texttt{makeQCHComp}$()$ in Algorithm \ref{algo:find_qch_comp} aids in this process. Notice that if any directed graph is passed to \texttt{makeQCHComp}$()$, it returns a collection of components. This is because in Step \ref{stp:highest_node}, a node $i$ with highest value of $\a_i$ is chosen. This definitely belongs to a QCH component as every node in itself is a QCH component. Then in Step \ref{stp:mpdp_path}, all nodes that have a path with quasi-concave MPDP parameters from this node is included into this component. Thus $(\tempnode,\temparc)$ at the end of a pass of the while loop is a DQCH by definition. The algorithm then repeats the process by considering the nodes which have not been visited in this process.

Now once such a collection of DQCH components is returned by \texttt{makeQCHComp}$()$ to \texttt{MACPartition}$()$ in Step \ref{stp:get_qch_comp}, Algorithm \ref{algo:mac_partition} checks every component and removes nodes that have out-neighbors not within that component (Step \ref{stp:remove_node_adqch}). Thus the nodes and arcs remaining (if any) at the end (Step \ref{stp:end_for_hq}) of each pass of the for-loop forms and A-DCQH component. The nodes and arcs added to $\node_\text{regroup}$ and $\node_\text{regroup}$ in Steps \ref{stp:v_regroup}, \ref{stp:a_regroup} have the property that they have outgoing arcs to some other components. Thus they are regrouped into NA-DQCH components in Step \ref{stp:regroup}. Note that any NA-DQCH component formed in this process cannot contain any A-DQCH because of the aforementioned property. Thus, the partition obtained at the end of Algorithm \ref{algo:mac_partition} satisfies all the requirements of Definition \ref{def:mac_partition} 
\bulletend
\end{remark}

Using such a partition, we can create a super graph of $\redflowx$ in
order to find a lower bound on the steady state social utility of the
dynamics. This is also useful as it reduces the number of variables
(significantly in some cases) to consider as we can group multiple
nodes into a super node and consider the super node as a whole rather
than considering the individual nodes that constitute it.

\begin{definition}\thmtitle{Maximal attractive component super graph or MAC-SG}
Suppose $\{\hill^{q} = (\snode^{q}, \sarc^{q})\}_{q \in [1,n]_\integer}$ is a MAC partition of $\flow = (\node,\arc)$. Construct a graph $\Gamma \ldef (\Lambda,\Pi)$ with the following property:
\begin{itemize}
	\item $\Lambda = \{q\}_{q \in [1,n]_\integer}$;
	\item $(q,r) \in \Pi$ if and only if $\exists$ $i \in \snode^{q}$ and $j \in \snode^{r}$ such that $(i,j) \in \arc$.
\end{itemize}
Such a graph $\Gamma$ is called a \emph{Maximal Attractive Component Super Graph} or \emph{MAC-SG} of the MAC partition $\{\hill^{q} = (\snode^{q}, \sarc^{q})\}_{p \in [1,n]_\integer}$ of $\flow$. The $q$'s are called super nodes.

The function $\texttt{MACSG}()$ takes in a MAC partition of $\flow$ and returns the corresponding MAC-SG.
\bulletend
\end{definition}

Let $\Gamma \ldef (\Lambda,\Pi)$ be such a MAC-SG of a MAC partition of $\redflowx$. 
Similar to the nodes, for a super node $q \in \Lambda$, we define the set of in-reachable super nodes of $q$ as $\inreach^{q} \ldef \{r \in \Lambda \,|\, \exists \text{ a directed path from } r \text{ to } q \text{ in } \Gamma\} \cup \{q\}$ and the set of out-reachable super nodes of $q$ as $\outreach^{q} \ldef \{r \in \Lambda \,|\, \exists \text{ a directed path from } q \text{ to } r \text{ in } \Gamma\} \cup \{q\}$.  We then let $\brho_{qr}$ denote the fraction of population moving from super node $q \in \Lambda$ to $r \in \Lambda$ and $\brho$ be the vector that accumulates all such $\brho_{qr}$ into a vector. Of course such a movement is only allowed between nodes if there is a directed path between them. Moreover, the total fraction of population that moves out of a super node cannot be more than the initial fraction that starts off in that node. Thus, 
\begin{subequations}
	\begin{align}
		\label{eq:min_opt_simp1}&\brho_{qr} \geq 0, \quad \forall r \in \outreach^{q}, \forall q \in \Lambda, \\
		\label{eq:min_opt_simp2}&\sum_{r \in \outreach^{q}} \brho_{qr} = \sum_{j \in \snode^{q}} \x^0_j, \quad \forall q \in \Lambda\,.
	\end{align}
\end{subequations}
Now, if we allow the entire population to fit in every super node, then we might end up with a conservative lower bound. To compute a more realistic lower bound, the fact that some nodes will not contain any population fraction in the steady state state must be taken into account. We formally define such nodes next and then provide a result to identify such nodes. 

\begin{definition}\label{def:eventually_empty}
\thmtitle{Eventually empty nodes}
Let the evolution of $\x$ be governed by a Nash convergent dynamics from an initial condition $\x^0$. Let $\Mcal{L}^+$ be the positive limit set of the trajectory. Then $i \in \node$ is said to be an eventually empty node if and only if $\bx_i = 0$, $\forall \, \bx \in \Mcal{L}^+$. \bulletend
\end{definition}

\begin{lemma}\label{lem:eventually_empty}
\thmtitle{Sufficient condition for being eventually empty for SSD, NBRD and NRPM}
Let $\x^0 \in \simplex^{|\node|}_\rho$ and let $\redflowx = (\rednodex,\redarcx)$ be an ICRG of $\graph = (\node, \edg)$. Let the evolution of $\x$ be governed by SSD, NBRD or NRPM. For a node $i \in \rednodex$, if $\exists \, j \in \neigh^i$  such that $(i,j) \in \redarcx$ and $(j,i) \notin \redarcx$, then $i$ is an eventually empty node.
\bulletend
\end{lemma}
A proof of this is provided in Appendix \ref{proof:eventually_empty}.
Thus, if a super node contains an eventually empty node, there exists
an inherent bound on the fraction of population that can stay in that
super node in the steady state state.  We discuss this in the following result.

\begin{lemma}\label{lem:super_node_bound_up}
\thmtitle{Upper bound on population fraction in super nodes for a Nash convergent dynamics}
Suppose $\x^0$ is the initial condition and $\redflowx$ is an ICRG of $\graph$. Let $\{\hill^q = (\snode^q,\sarc^q)\}_{q \in [1,n]_\integer}$ be a MAC partition of $\redflowx$ and let $\Gamma = (\Lambda,\Pi)$ be the corresponding MAC-SG. Let $\poslim$ be the positive limit set of the trajectories of a Nash convergent dynamics. For a super node $q \in \Lambda$, let the set of eventually empty nodes in $q$ be denoted by $\M^{q} \ldef \{i \in \snode^{q} \,|\, i \text{ is eventually empty}\}$. Let
\begin{align}
	\a^\text{max}_{q} \ldef \max_{i \in \M^{q}} \a_i, \quad \text{ if } \M^{q} \neq \varnothing\,,
\end{align}
and finally, let
\begin{equation}
	\rhom_{q} \ldef 
	\begin{cases}
		0, & \text{if } \M^{q} = \snode^{q}\\
		\displaystyle \sum_{i \in \snode^{q} \setminus \M^{q}} \lev^{-1}_i\left( \a^\text{max}_{q} \right), & \text{if } \snode^{q} \supset \M^{q} \neq \varnothing,\\
		1, &  \text{if } \M^{q} = \varnothing\,.
	\end{cases}
	\label{eq:rho_max}
\end{equation}
Then, $\forall \, \bx \in \poslim$
\begin{equation}
	\sum_{i \in \snode^q} \bx_i \leq \rhom_q, \,\, \forall \, q \in \Lambda\,.
	\label{eq:snode_up_bound}
\end{equation}
\bulletend
\end{lemma}
We present a proof in Appendix \ref{proof:super_node_bound_up}.

Recall that $\optu_{\dumnode}(\rho)$ is the solution of the optimization problem $\prob_1(\dumnode,\rho)$ in \eqref{eq:P_no_restrict}. The optimum of the following optimization problem then provides a lower bound to the steady state social utility.

\begin{equation}
	\begin{split}
		\prob_5: \qquad  & 
		\begin{split}		
			U_\text{min} \ldef \min_{\xi,\brho} & \quad \sum_{q \in \Lambda} \optu_{\snode^{q}} \left( \xi_{q} \right)\\
			\text{s.t.} & \quad \xi_{q} = \sum_{r \in \inreach^{q} } \brho_{rq}, \forall q \in \Lambda,\\
			& \quad \eqref{eq:min_opt_simp1}, \eqref{eq:min_opt_simp2}, \\
			& \quad \brho_{q} \leq \rhom_{q}, \forall q \in \Lambda \,.
		\end{split}
	\end{split}
\end{equation}

The following lemma, proof of which is given in Appendix \ref{proof:p1_concave_rho} shows concavity of the cost function in $\prob_5$.
\begin{lemma} \label{lem:p1_concave_rho}
Consider a fixed $\dumnode$ and the function $\optu_{\dumnode}(\rho)$ which is the optimum of $\prob_1(\dumnode,\rho)$ in \eqref{eq:P_no_restrict}. Then $\optu_{\dumnode}(\rho)$ is concave in $\rho \in \real_+$.
\bulletend
\end{lemma}
The constraints set of $\prob_5$ is a polyhedron. Thus it is a well
known fact that the global optimum will occur at an extreme point.
Any standard method such as \cite{MED-LGP:1977:all_ext_pt,
  THM-DSR:1980:ext_pt_survey, HPB:1985:cancave_min_finite} can be used
to solve this problem.

The main sequence of steps, required to compute a lower bound on the steady state social utility, described in this section is listed in Algorithm \ref{algo:u_min} and the result regarding the same is stated in the next to conclude the section. 

\begin{algorithm}
\caption{Compute lower bound for steady state social utility}
\label{algo:u_min}
\SetAlgoLined
\KwData{$\x^0, \flow$}
\KwResult{$\umin$}


$\redflowx \gets \texttt{reduceGraph}(\flow,\x^0)$ \label{stp:final_algo_red_graph}\\
$\hill \gets \texttt{MACPartition}(\redflow_0)$\label{stp:final_algo_macp}\\
\nonl \mycomment{individual components are $\hill^q = (\snode^q, \sarc^q)$}\\
$\Gamma \gets \texttt{MACSG}(\hill)$  \mycomment{$\Gamma = (\Lambda,\Pi)$}\label{stp:final_algo_macsg}\\
$\umin \gets$ solution of $\prob_5$ \label{stp:u_min}\\
\end{algorithm}

\begin{theorem} \label{th:lower_bound}
\thmtitle{Lower bound on steady state social utility for SSD, NBRD and NRPM}
Let $\x^0 \in \simplex^{|\node|}_\rho$. Then
\begin{align*}
	\umin \leq \bussd; \,\, \umin \leq \bunbrd; \,\, \umin \leq \bunrpm\,.
\end{align*}
\bulletend
\end{theorem}
A proof of this is presented in Appendix \ref{proof:lower_bound}

\section{Conclusion}

We provided sufficient conditions on the network under which there
exists a unique Nash equilibrium for the stratified population
game. We also provided algorithms to reduce the graph using the
initial condition without affecting the population evolution. We then
provided algorithms to partition the reduced graph and utilized the
conditions for unique Nash equilibrium to provide upper and lower
bounds on the steady state social utility for SSD, NBRD and NRPM.

Future work includes utilizing the dynamics to further reduce the
graph and refine the bounds on the steady state social utility. We
would also like to extend the ideas further to compute conditions
under which the bounds coincide and utilize this knowledge to compute
the steady state state.

\bibliographystyle{IEEEtran}
\bibliography{myref}

\begin{IEEEbiography}[{\includegraphics[width=1in,height=1.25in,clip,keepaspectratio]{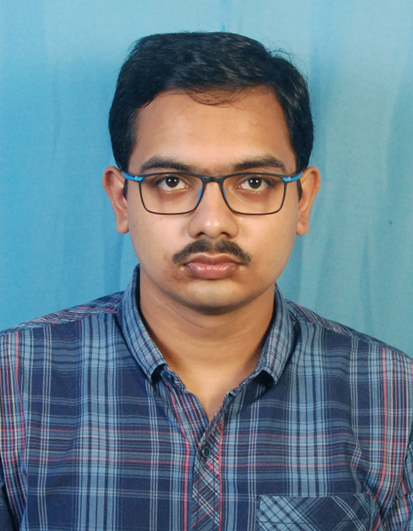}}]{Nirabhra Mandal} received the B.Tech. degree in Electrical Engineering from Institute of Engineering and Management, Salt Lake, Kolkata, India in 2017. Since 2018, he is pursuing the M.Tech(Res) degree from the Department of Electrical Engineering at the Indian Institute of Science. His research interests include multi-agent systems, population games, evolutionary dynamics on networks and non-linear control.
\end{IEEEbiography}
\vskip 0pt plus -1fil
\begin{IEEEbiography}[{\includegraphics[width=1in,
    height=1.25in,clip,keepaspectratio] {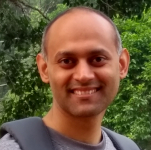}}]{Pavankumar
    Tallapragada} (S'12-M'14) received the B.E. degree in
  Instrumentation Engineering from SGGS Institute of Engineering $\&$
  Technology, Nanded, India in 2005, M.Sc. (Engg.) degree in
  Instrumentation from the Indian Institute of Science in 2007 and the
  Ph.D. degree in Mechanical Engineering from the University of
  Maryland, College Park in 2013. He was a Postdoctoral Scholar in the
  Department of Mechanical and Aerospace Engineering at the University
  of California, San Diego from 2014 to 2017. He is currently an
  Assistant Professor in the Department of Electrical Engineering and
  the Robert Bosch Centre for Cyber Physical Systems at the Indian
  Institute of Science. His research interests include networked
  control systems, distributed control, multi-agent systems and
  networked transportation systems.
\end{IEEEbiography} 

\appendices

\section{Simulations and Analysis} \label{sec:sims}

In this section, we highlight some properties of SSD, NBRD and NRPM
and illustrate some results stated in Parts I and II using some
simulations. We used CVXPY \cite{SD-SB:2016:cvxpy,
  AA-RV-SD-SB:2018:cvxpy_rewriting} for solving the optimization
problems. We performed all simulations in a python3 programming
language environment on a standard laptop with intel 10th generation
Core i5 processor.

\subsection{Myopic Coordination can be Worse than Myopic %
  Selfish Behavior}

Even though SSD, NBRD and NRPM converge asymptotically to some point
in the set of Nash equilibria, they may not converge to the same state
in general. Moreover, even though SSD, NBRD and NRPM display
increasing levels of coordination among the agents, it is not always
true that NBRD performs better than SSD nor that NRPM performs better
than NBRD. This is because even though there is coordination within a
node in NBRD and coordination of the entire population in NRPM, they
are still essentially network restricted gradient ascent dynamics.

In this section we illustrate two situations where the act
of coordination does not result in better social utility. For this
section we assume the cumulative payoff functions to be of the form
$\pe_i(\x_i) \ldef -0.5\x^2_i - \a_i\x_i$. This is the uniform water
tank model presented in \cite{NM-PT:2020:ifac}.

\subsubsection{SSD outperforms NBRD}

\begin{figure}
  \begin{center}
    \begin{tabular}{cc}
      \includegraphics[scale=0.25]{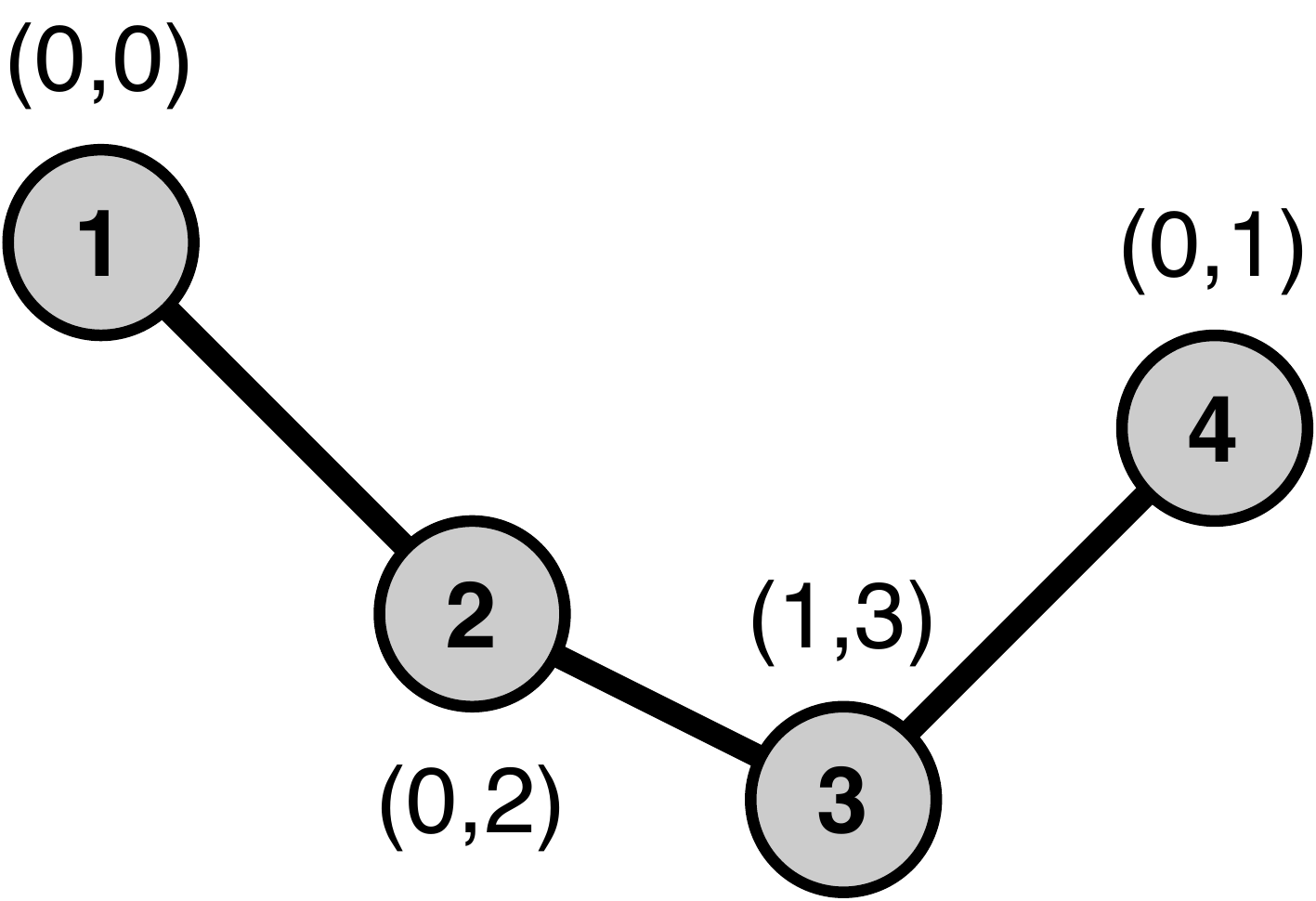} %
      & \includegraphics[scale=0.25]{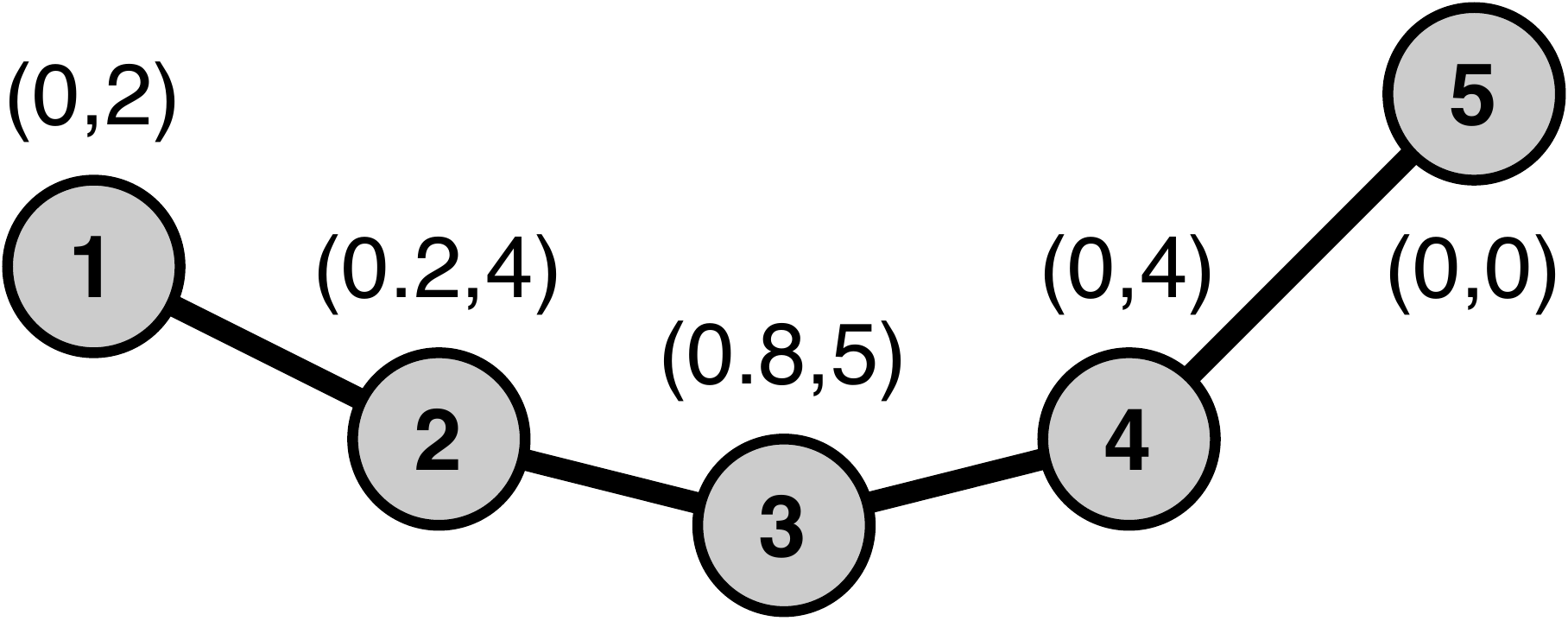} \\
      \hypertarget{subfig:ssd_beats_nbrd}{(a)} %
      & \hypertarget{subfig:nbrd_beats_nrpm}{(b)}  
    \end{tabular}
  \end{center}
  \caption{Graph structure, initial state and MPDP's for showing that
    coordination is not always good. The tuple $(.,.)$ near each node
    $i$ represents $(\x^0_i,\a_i)$. (a) SSD outperforms NBRD. (b) NBRD
    outperforms NRPM.}
  \label{fig:coop_anomaly}
\end{figure}

\begin{figure*}
  \begin{center}
    \begin{tabular}{ccc}
      \includegraphics[scale=0.35]{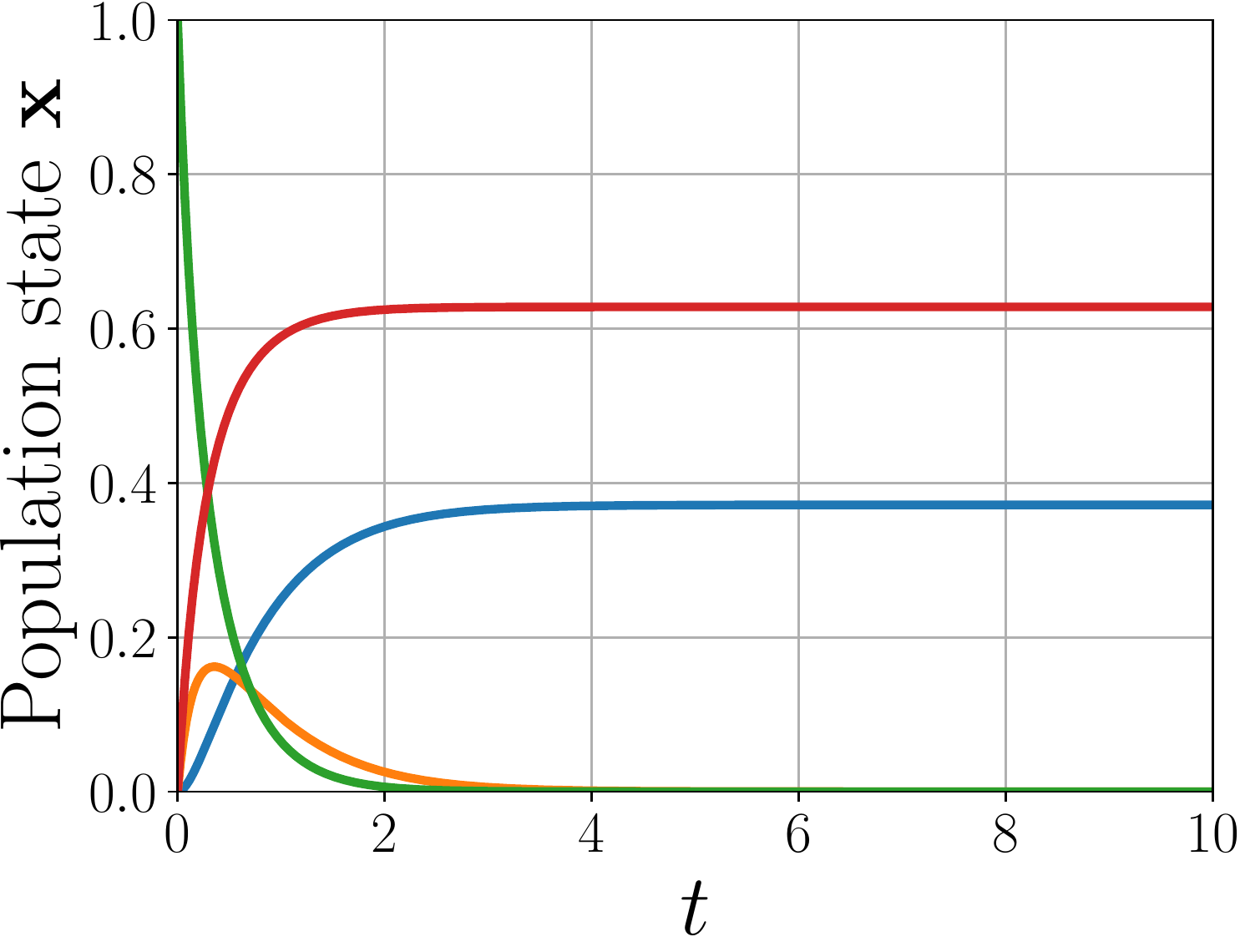} %
      & \includegraphics[scale=0.35]{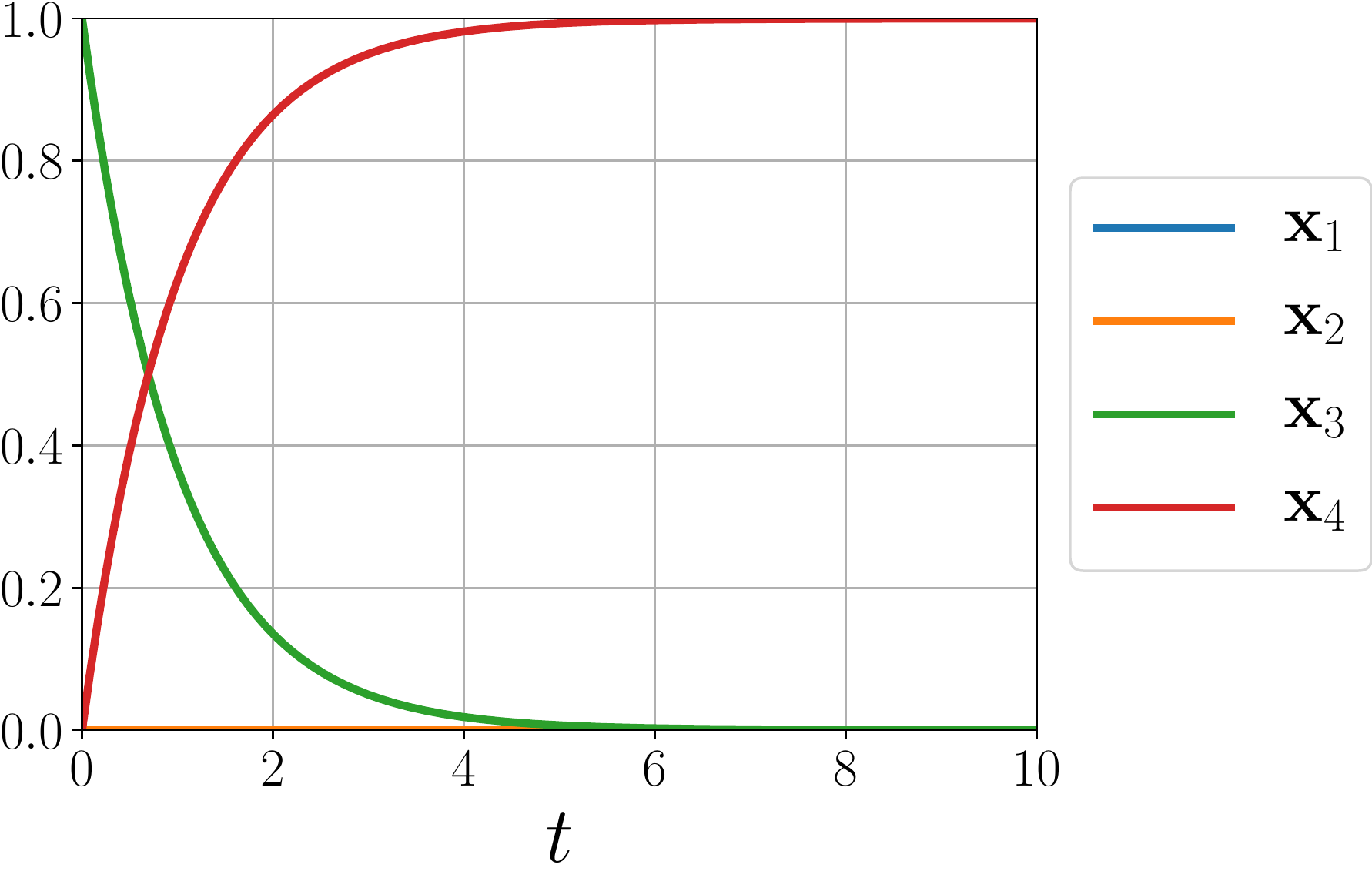} %
      & \includegraphics[scale=0.35]{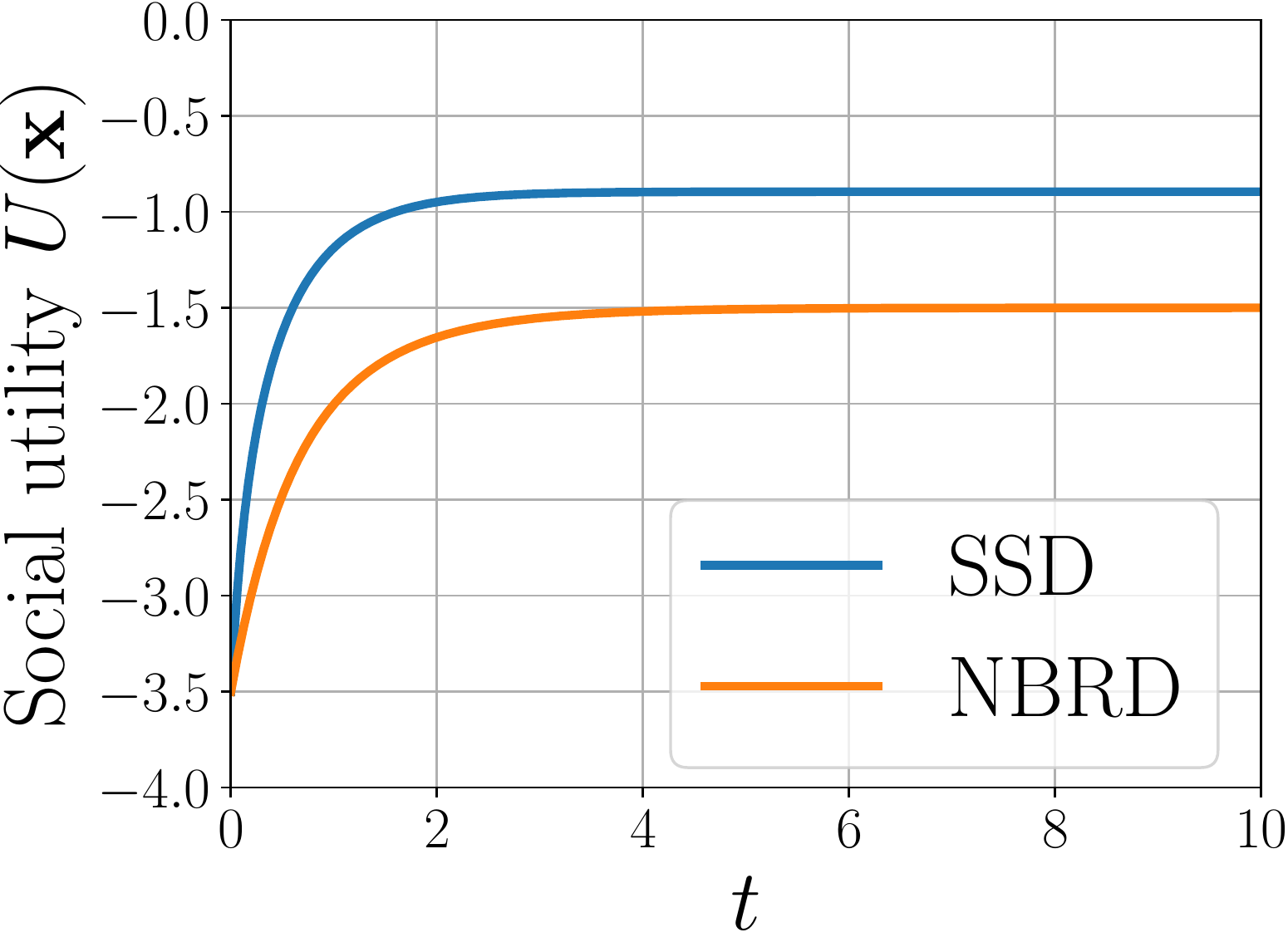}\\
      (a) & (b) & (c)
    \end{tabular}
  \end{center}
  \caption{SSD performs better than NBRD. The first two plots share a
    common label and a common legend. (a) Evolution of population
    state under SSD. (b) Evolution of population state under NBRD. (c)
    Evolution of social utility.}
  \label{fig:ssd_better_than_nbrd}
\end{figure*}

To demonstrate such a phenomenon we consider a four node graph with
initial state and MPDP's described in Figure
\ref{fig:coop_anomaly}\hyperlink{subfig:ssd_beats_nbrd}{(a)}. Note
that in the case of NBRD, the outflow from node 3 is always directed
only to node 4. Hence the population misses out on a better payoff in
node 1. For SSD, on the other hand, there are outflow from node 3 both
towards nodes 2 and 4. Thus a fraction of the population gets a chance
to move to node 1. Hence the population as a whole receives a higher
utility. The simulation results are shown in Figure
\ref{fig:ssd_better_than_nbrd}.

\subsubsection{NBRD outperforms NRPM}

\begin{figure*}
  \begin{center}
    \begin{tabular}{ccc}
      \includegraphics[scale=0.35]{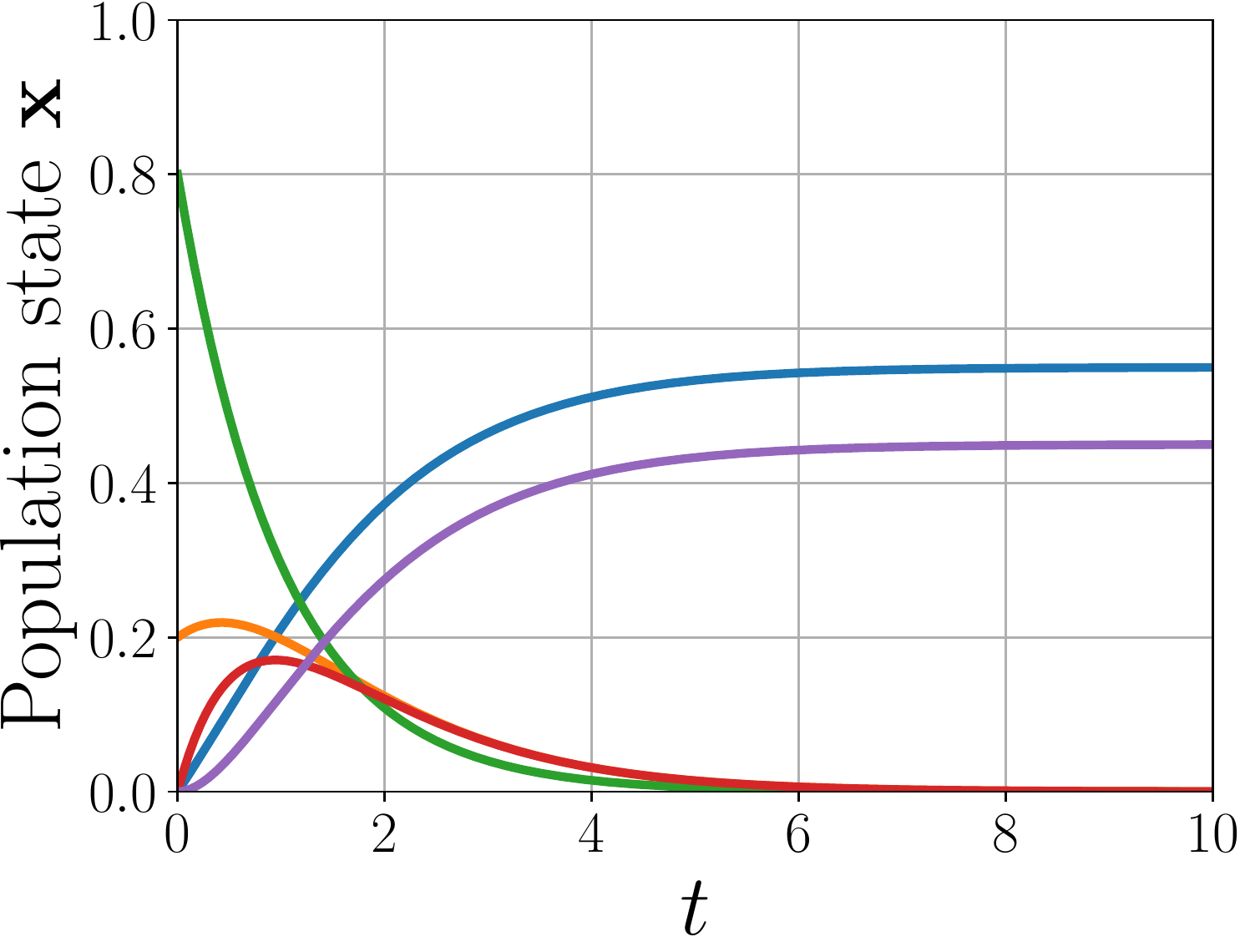} & \includegraphics[scale=0.35]{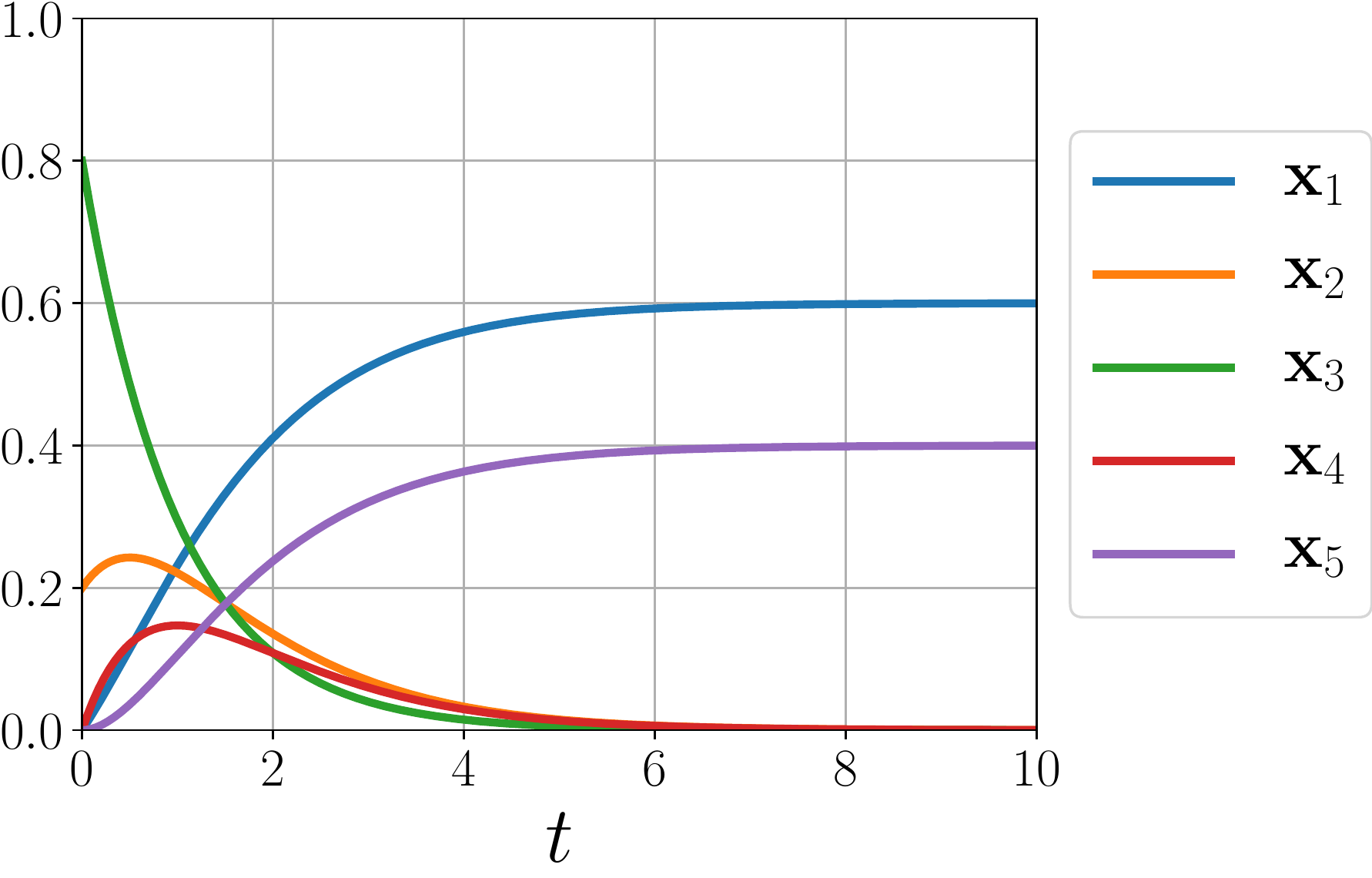} & \includegraphics[scale=0.35]{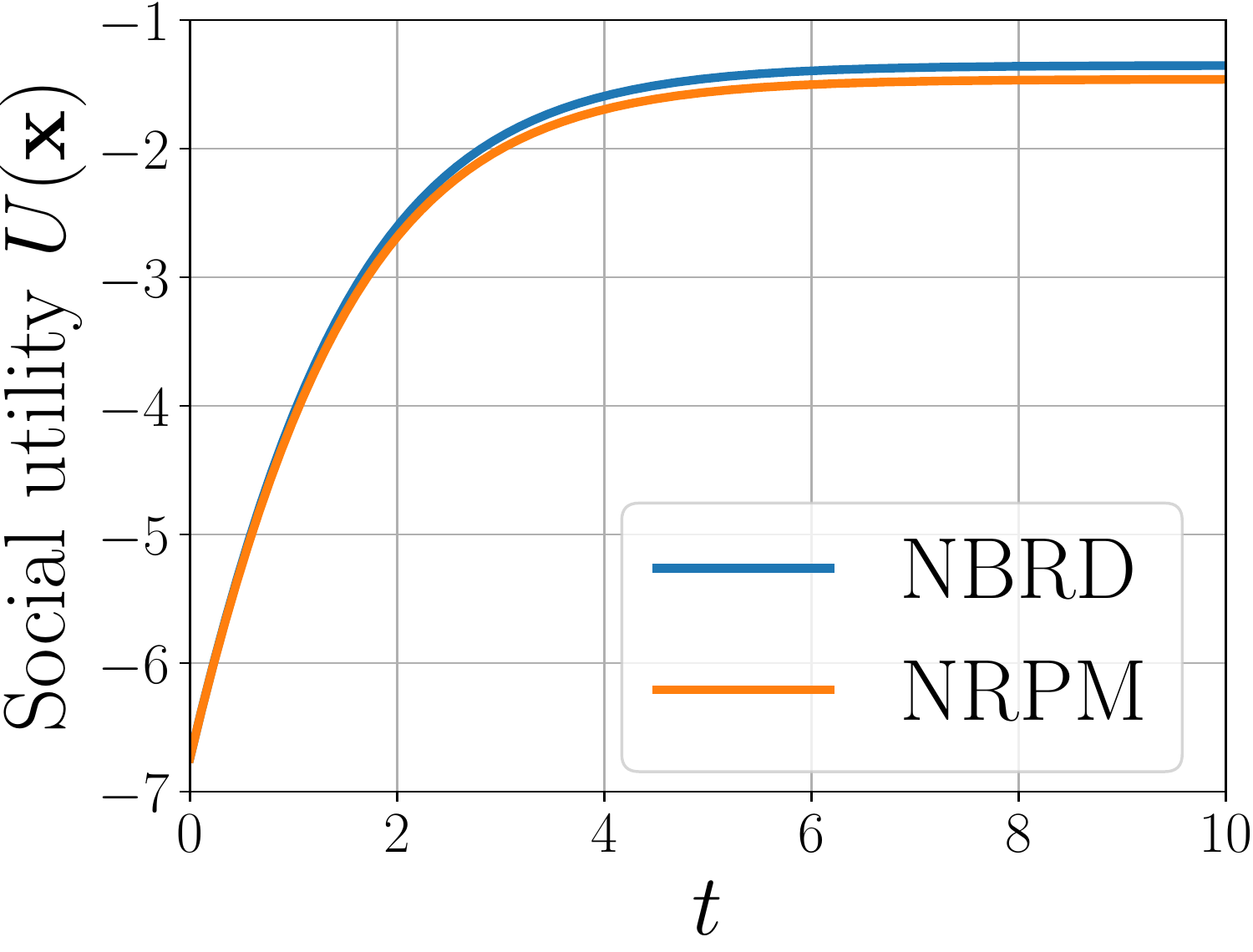}\\
      (a) & (b) & (c)
    \end{tabular}
  \end{center}
  \caption{NBRD performs better than NRPM. The first two plots share a
    common label and a common legend. (a) Evolution of population
    state under NBRD. (b) Evolution of population state under
    NRPM. (c) Evolution of social utility.}
  \label{fig:nbrd_better_than_nrpm}
\end{figure*}

To demonstrate this phenomenon we consider a five node graph with
initial state and MPDP's described in Figure
\ref{fig:coop_anomaly}\hyperlink{subfig:nbrd_beats_nrpm}{(b)}. Note
that in the case of NBRD, the population fraction in node 3 is unaware
that fraction in node 2 is also simultaneously revising is revising
its choices. Thus, there is a greater outflow from node 3 to node 4
than to node 3. However, in the case of NRPM, the fraction in node 3
is aware of the fact that the fraction in node 2 will move to node
1. Hence there are equal outflows from node 3 to nodes 2 and 4. Thus a
larger fraction of population gets to visit node 5 in case of NBRD
than for NRPM. This causes the steady state social utility of NBRD to be
greater than NRPM. Simulation results illustrating this are provided
in Figure \ref{fig:nbrd_better_than_nrpm}.

\subsection{Bounds on Steady-State Social Utility}

\begin{figure*}
	\begin{center}
	\begin{tabular}{ccc}
		 \includegraphics[scale=0.24]{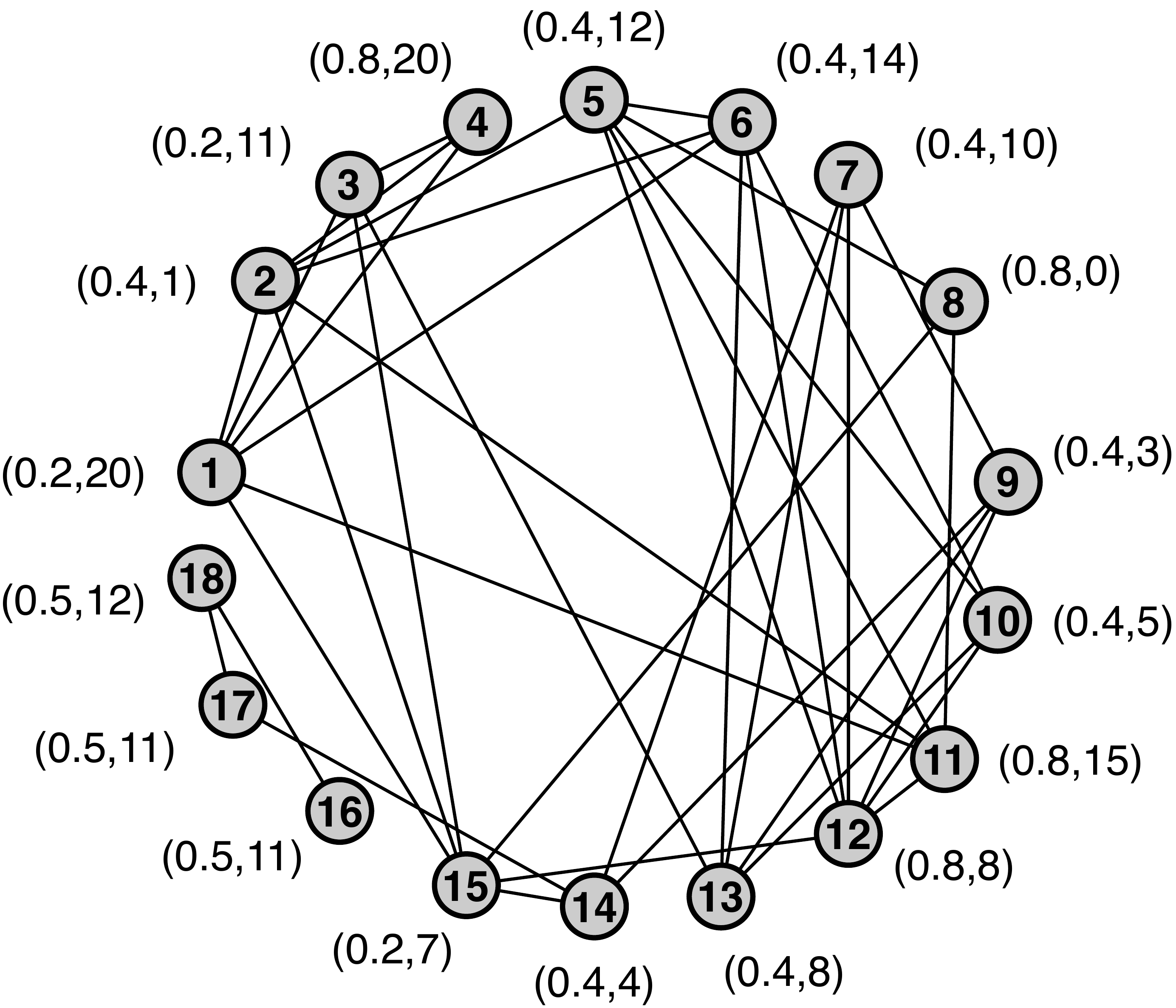} &\includegraphics[scale=0.24]{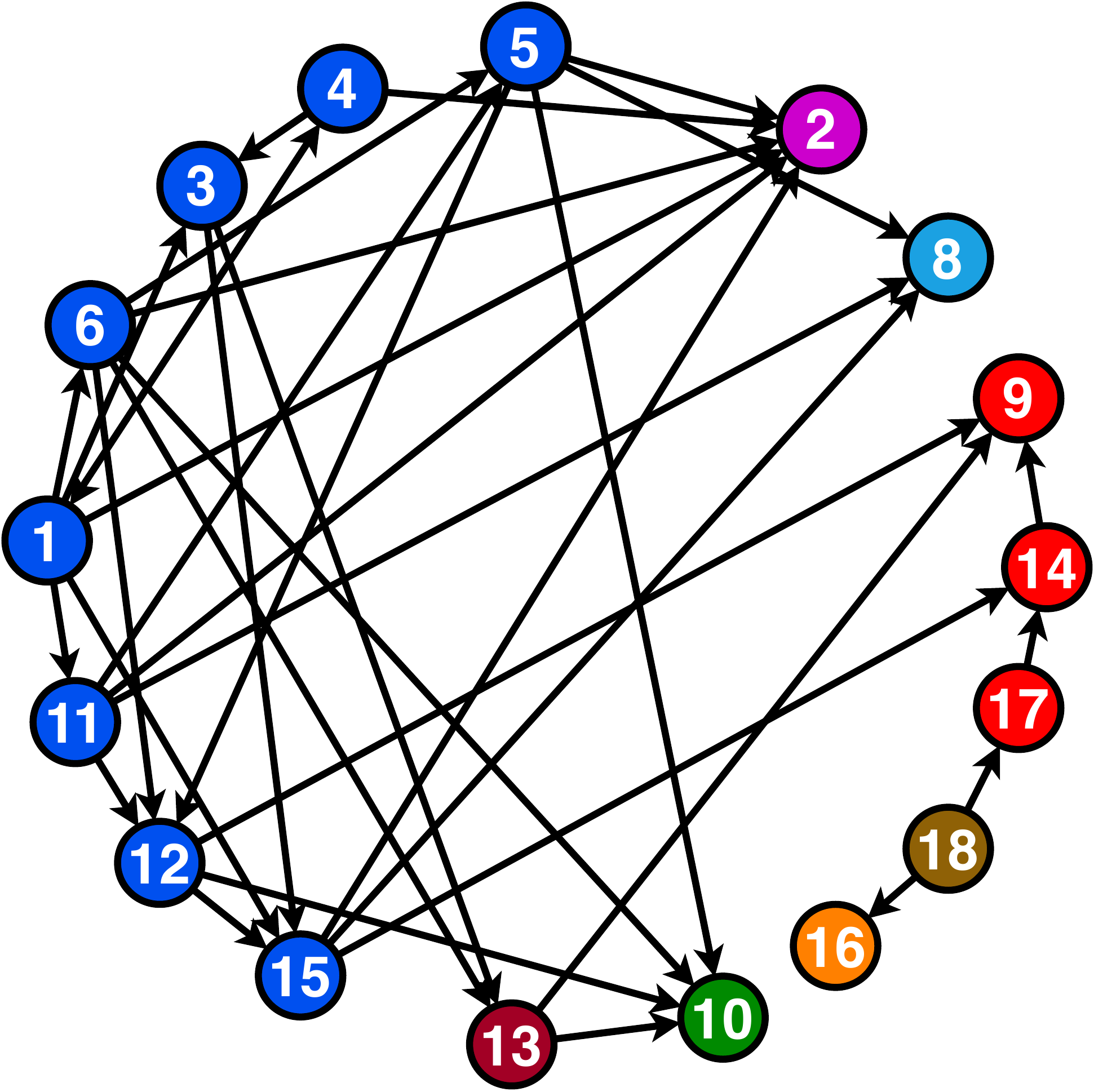} & \includegraphics[scale=0.34]{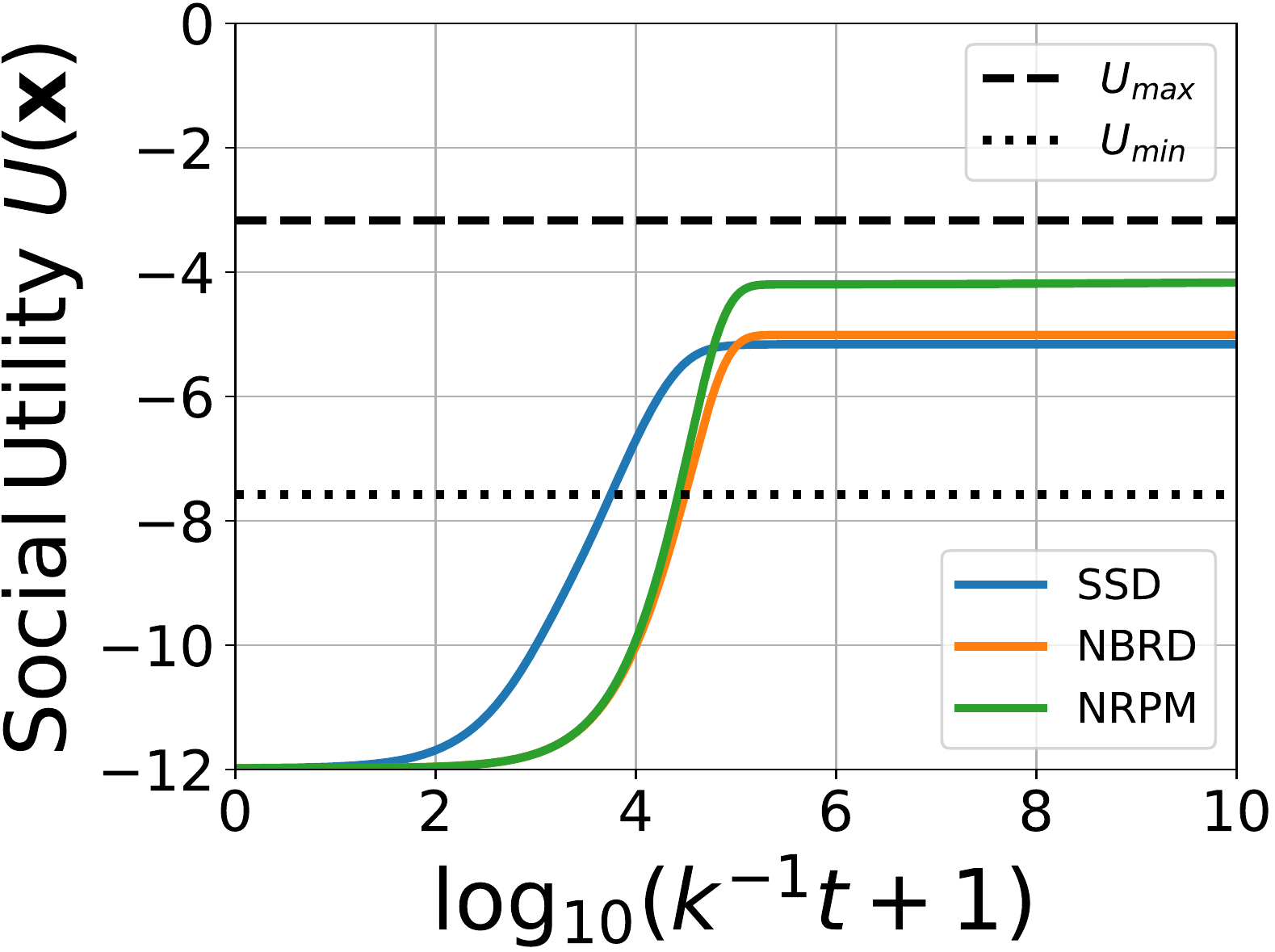}\\
		\hypertarget{subfig:18_node_graph}{(a)} & \hypertarget{subfig:18_node_partition}{(b)} & \hypertarget{subfig:18_node_rest1}{(c)}
	\end{tabular}
	
	\begin{tabular}{ccc}
		\includegraphics[scale=0.34]{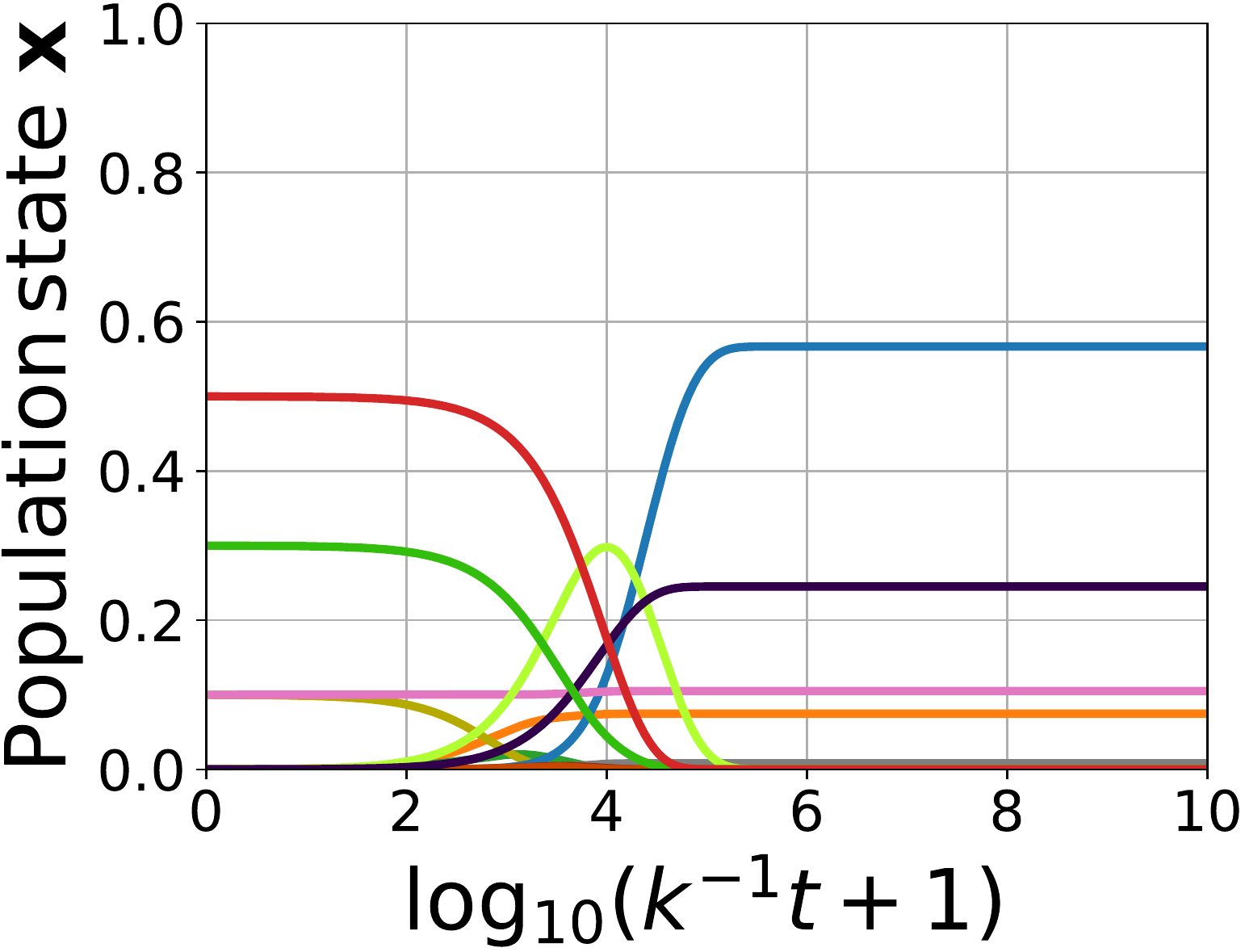} & \includegraphics[scale=0.34]{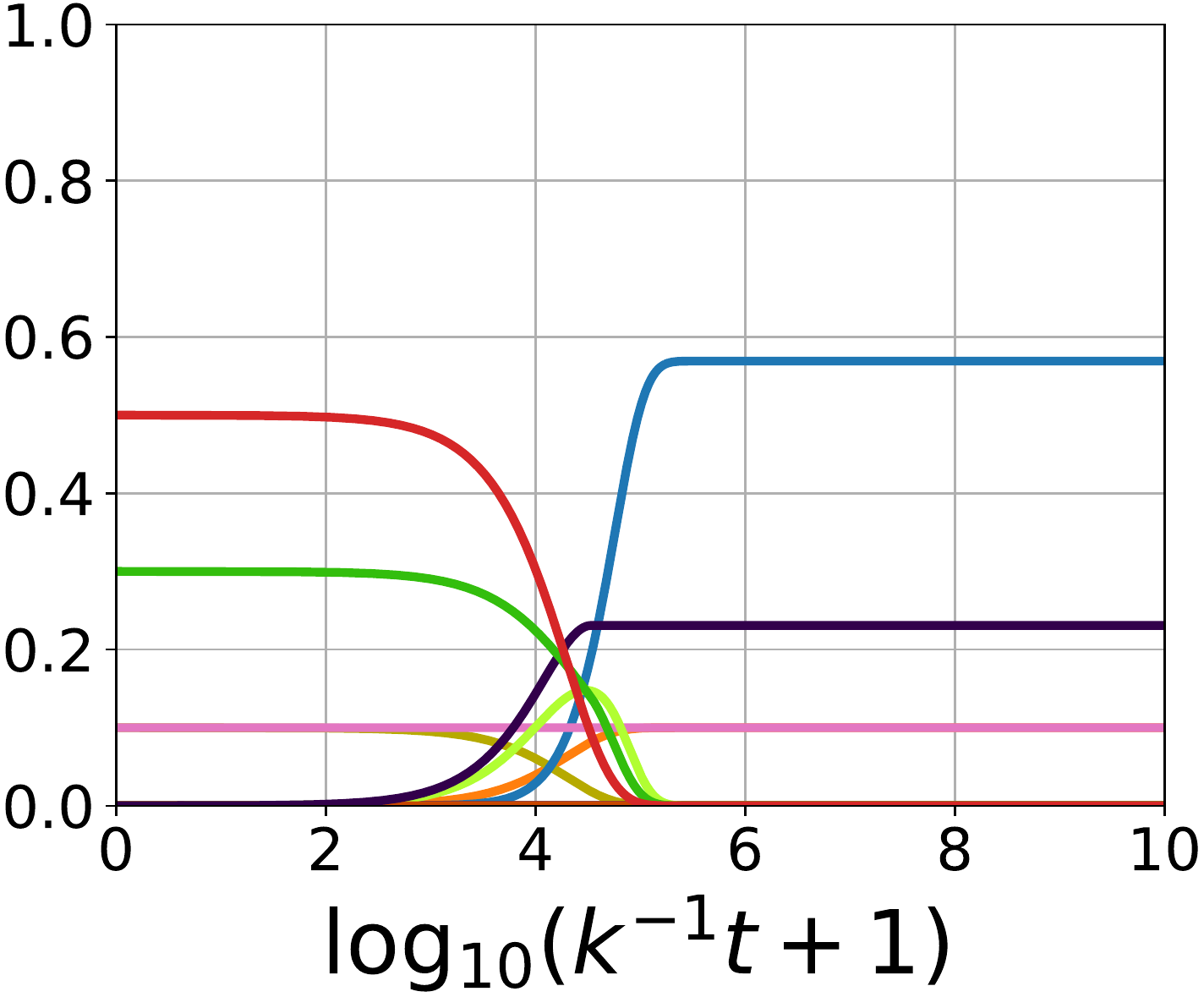} & \includegraphics[scale=0.34]{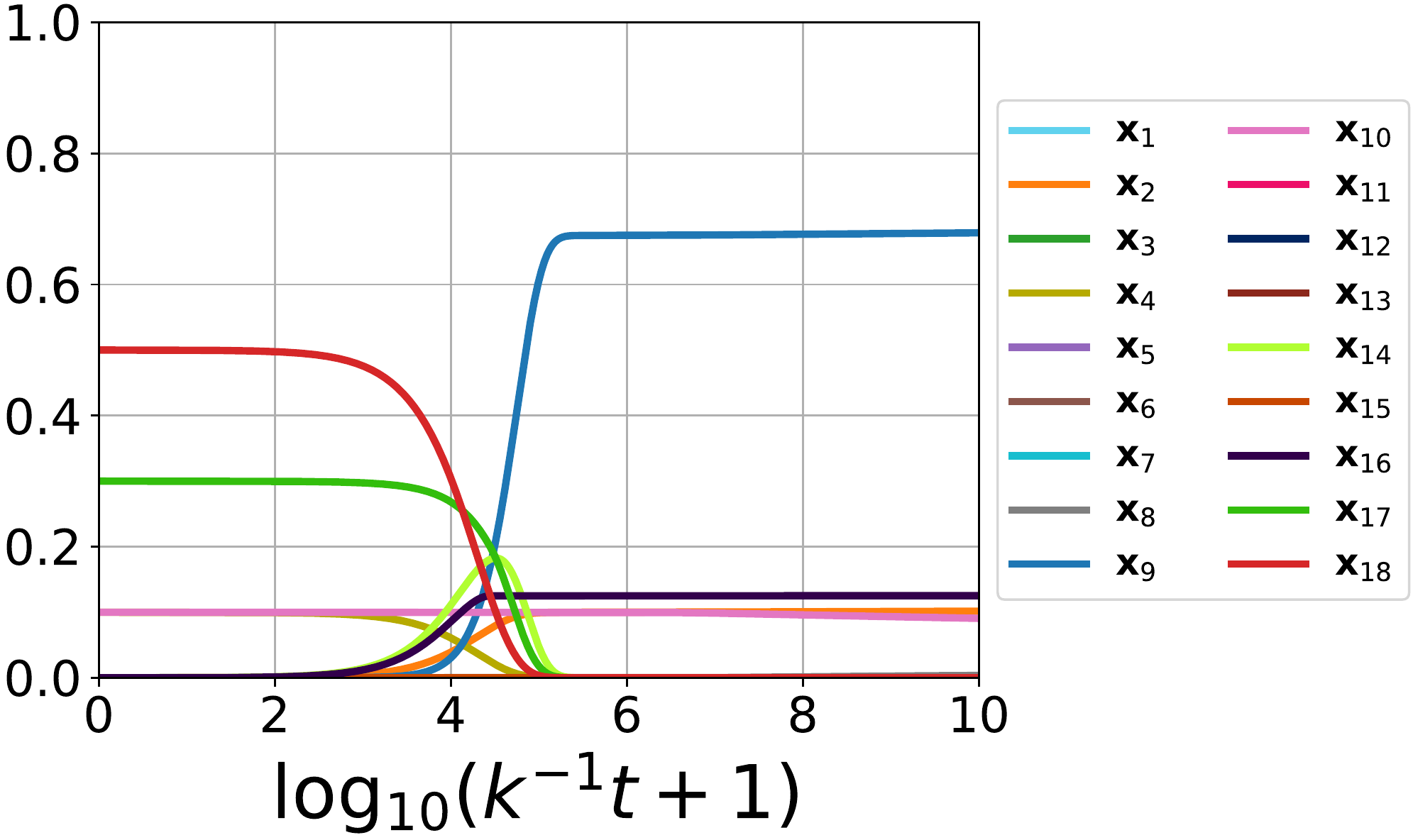}\\
		(d) & (e) & \hypertarget{subfig:18_node_rest2}{(f)}
	\end{tabular}
	\end{center}
	\caption{Simulation with $\pe_i(\x_i) \ldef -\alpha_i \,\x_i^2 - \beta_i\,\x_i$. The last three plots share a common label and a common legend. Here $k \approx 5 \times 10^{-5}$. (a) Graph structure with 18 nodes. The tuple (.,.) around each node $i$ represents $(\alpha_i,\beta_i)$. (b) Corresponding ICRG. Nodes in same D-QCH components have same color and have been grouped together. Nodes in different D-QCH components are colored differently.  (c) Evolution of social utility and bounds on steady-state social utility. (d) Evolution of population state under SSD. (e) Evolution of population state under NBRD. (f) Evolution of population state under NRPM.}
    \label{fig:big_sim}
\end{figure*}

Here we consider an 18 node graph and quadratic cumulative payoff
functions and verify the upper and lower bounds on the steady-state
social utility. The graph structure $\graph$ and cumulative payoff
function details are provided in Figure
\ref{fig:big_sim}\hyperlink{susubfig:18_node_graph}{(a)}. Initial
population is distributed between nodes $4, 10, 17$ and $18$ as
$\x_4^0 = \x_{10}^0 = 0.1$, $\x^0_{17} = 0.3$, $\x^0_{18} = 0.5$ and
$\x_i^0 = 0$,
$\forall i \in [1,18]_\integer \setminus \{4,10,17,18\}$.  The
corresponding ICRG $\redflowx$ is provided in Figure
\ref{fig:big_sim}\hyperlink{susubfig:18_node_partition}{(b)}. Note
that $\redflowx$ has only $17$ nodes and not $18$ (node $7$ has been
removed). The only bi-directional arc in $\redflowx$ is between nodes
$1$ and $4$. Rest of the arcs are all uni-directional. Thus there is
over $50\%$ reduction of $\delta_{ij}$ variables from $\graph$ to
$\redflowx$. The A-DQCH components are made up of the node sets
$\{2\}, \{8\}, \{9,14,17\}, \{10\}, \{16\}$ and corresponding
arcs. The NA-DQCH components are made up of the node sets
$\{1,3,4,5,6,11,12,15\}, \{13\}, \{18\}$ and corresponding arcs. To
compute the upper and lower bounds the total time taken was
$\approx 3.48$ sec.

The simulation results are provided in Figures
\ref{fig:big_sim}\hyperlink{susubfig:18_node_rest1}{(c)} -
\ref{fig:big_sim}\hyperlink{susubfig:18_node_rest2}{(f)}. The time
axis in each of these diagrams is represented as a log scale,
\emph{i.e.} each tick on the horizontal axis represents
$\log(k^{-1}t + 1)$ rather than $t$ (the value of $k$ is given in
Figure \ref{fig:big_sim}). This is done to magnify the initial time
frame where the main redistribution occurs and shrink the later time
frame. The average time taken to solve the optimization problems for
NBRD and NRPM were $\approx 0.32$ sec and $\approx 0.16$ sec
respectively. The total time required to complete the simulations for
SSD, NBRD and NRPM were $\approx 108.43$ sec,
$\approx 3.19 \times 10^{4}$ sec and $\approx 1.6 \times 10^{4}$ sec
respectively. The upper and lower bounds obtained were $-3.17$ and
$-7.58$ respectively while the actual steady-state social utility for
SSD, NBRD and NRPM were $-5.16$, $-5.01$ and $-4.17$ respectively.

\section{}
\subsection{Proof of Results on Unique Nash Equilibrium}

\subsubsection{Proof of Lemma \ref{lem:monotone_path}} \label{proof:monotone_path}
\begin{proofa}\emph{(By contradiction)} Suppose $\exists \, k \in [1,n]_\integer$, $p \in [1,k-1]_\integer$ and $q \in [k+1,n]_\integer$ such that $a_{\pi(p)} > \a_{\pi(p+1)}$ and $\a_{\pi(q-1)} < \a_{\pi(q)}$. But $\pi(p),\pi(q) \in \gpath(i,j)$ which is a path of quasi-concave MPDP's and one of the previous inequalities violates the condition that $\forall \, r \in [p,q]_\integer$, $\a_{\pi(r)} \geq \min \{\a_{\pi(p)},\a_{\pi(q)}\}$. This is a contradiction and hence the assumption was incorrect.
\end{proofa}

\subsubsection{Proof of Lemma \ref{lem:same_lev}} \label{proof:same_lev}
\begin{proofa}
  Choose any $i,j \in \supp(\x)$ and let $\gpath(i,j)$ be a path
  between $i$ and $j$ with quasi-concave MPDP's. Existence of such a
  path is guaranteed as the graph $\graph$ is a QCH. Let
  $\pi(k) \ldef \gpath_k(i,j)$. As $\pi(k) \in \gpath(i,j)$, by Lemma
  \ref{lem:monotone_path}
  $\a_i \leq \a_{\pi(2)} \leq \cdots \a_{\pi(k)}$ or
  $\a_{\pi(k)} \geq \cdots \geq \a_{\pi(n-1)} \geq \a_j$. Without loss
  of generality, suppose that
  $\a_i \leq \a_{\pi(2)} \leq \cdots \a_{\pi(k)}$.

  Now, consider the nodes $\pi(1) = i$ and $\pi(2)$. As $\x_i > 0$ and
  $\a_i \leq \a_{\pi(2)}$, we have
  \begin{equation*}
    \lev_{\pi(2)}(0) \geq \lev_i(0) > \lev_i(\x_i) \geq
    \lev_{\pi(2)}(\x_{\pi(2)}) ,  
  \end{equation*}
  where we have used the fact that $\lev_k(.)$ is a strictly
  decreasing function $\forall \, k$, which implies
  $\lev_k(0) > \lev_k(y)$ iff $y > 0$. Further, the last inequality
  comes from the fact that $\x \in \ne^{|\node|}_\rho$. Thus
  $\x \in \ne^{|\node|}_\rho$, $i \in \supp(\x)$ and
  $\a_i \leq \a_{\pi(2)}$ together implies $\x_{\pi(2)} > 0$. Now,
  repeating this argument for every pair of nodes $(\pi(r), \pi(r+1))$
  for $r \in \{ 1, \ldots k-1\}$, we conclude that $\x_{\pi(k)} > 0$,
  that is $\pi(k) \in \supp(\x)$.

  Thus, $\forall \, \pi(k) \in \gpath(i,j)$, $\pi(k) \in
  \supp(\x)$. Then as $\x \in \ne^{|\node|}_\rho$,
  $\lev_{\pi(k)}(\x_{\pi(k)}) = \lev_{\pi(l)}(\x_{\pi(l)})$
  $\forall \, \pi(k),\pi(l) \in \gpath(i,j)$. Thus
  $\lev_i(\x_i) = \lev_j(\x_j)$. As $i$ and $j$ were chosen
  arbitrarily, the proof of the lemma is complete.
\end{proofa}

\subsubsection{Proof of Theorem \ref{th:unique_nash}}\label{proof:unique_nash}

\begin{proofa}
  First note that $U(\x)$ is a strictly concave function in $\x$ and
  $\prob_1$ is always feasible as $\simplex^{|\node|}_\rho$ is
  non-empty. Also $\prob_1$ is a strictly convex program and hence has
  a unique optimizer \cite{SB-LV:2004:cvx}. Thus, it suffices to show
  that if $\x \in \ne^{|\node|}_\rho$ then $\x$ also optimizes
  $\prob_1$.

  If $\rho = 0$ then the result is trivially true. So, now suppose
  $\rho > 0$. The Lagrangian for $\prob_1$ can be written as
\begin{equation*}
  \lag_1 = \sum_{i \in \node} \pe_i(\x_i) - \lambda \left(\sum_{i \in \node} \x_i - \rho \right) + \sum_{i \in \node} \mu_i \x_i,
\end{equation*}
where $\lambda$ and $\{\mu_i \geq 0\}_{i \in \node}$ are the Lagrange
multipliers.  The KKT conditions for $\prob_1$ include
\begin{subequations}\label{eq:kkt_P1}
  \begin{align}
    \label{eq:kkt_grad}
    &\lev_i(\x_i) - \lambda + \mu_i = 0, \, \forall i \in \node,\\
    \label{eq:kkt_comp_slack}
    &\mu_i \x_i = 0, \, \forall i \in \node \,.
  \end{align}
\end{subequations}
Now, let $\x \in \ne^{|\node|}_\rho$. Since the graph is a QCH, we
know from Lemma~\ref{lem:same_lev} that
\begin{equation}
  \lev_i(\x_i) = H, \quad \forall \, i \in
  \supp(\x), \label{eq:QCH_lev}
\end{equation}
for some $H$. Clearly, if $\rho > 0$, problem $\prob_1$ satisfies
Slater's condition. Thus, we will show that
$\x \in \ne^{|\node|}_\rho$ is the unique optimizer of $\prob_1$ by
showing that for $\x$ there exist Lagrange multipliers that
satisfy~\eqref{eq:kkt_P1} and the feasibility constraints. Note that
as $\x \in \ne^{|\node|}_\rho$, $\x$ is a feasible solution for
$\prob_1$. Now, we set
\begin{align*}
  \lambda^* = H, \quad %
  &\mu^*_i = 0, \ \forall i \in \supp(\x), \\
  &\mu^*_j = H - \lev_j(0), \ \forall j \notin \supp(\x).
\end{align*}
As $\x \in \ne^{|\node|}_\rho$, we can say from~\eqref{eq:QCH_lev}
that $\mu^*_j \geq 0$, $\forall j \notin \supp(\x)$. Thus for each
$i \in \node$, $\mu^*_i \geq 0$ and satisfies
\eqref{eq:kkt_comp_slack}. Also, we can directly verify that
$(\x, \lambda^*, \mu^*)$ satisfy~\eqref{eq:kkt_grad}. Thus $\x$ must
be the unique optimizer of $\prob_1$. This proves the theorem.

\end{proofa}

\subsection{Proof of Results on Bounds on Social Utility}

\subsubsection{Proof of Lemma \ref{lem:arc_remove}} \label{proof:arc_remove}
\begin{proofa}
Since $\lev_i(.)$ is strictly decreasing $\forall i \in \node$, we
  can say that $\forall (i,j) \in \arc$ such that
  $\lev_i(\theta_i) \geq \lev_j(0)$
  \begin{equation*}
    \lev_i(\x_i) \geq \lev_i(\theta_i) \geq \lev_j(0) \geq
    \lev_j(\x_j), \quad \forall \x_i \in [0, \theta_i] , \ \forall
    \x_j \geq 0 . 
  \end{equation*}
  Now, for SSD and NBRD, \eqref{eq:delta_zero_cond} follows from strong positive correlation of
  SSD and NBRD with $\lev(.)$ from Theorems \ref{th:ssd_conv} and \ref{th:NBRD_full}.
  
  We prove the condition for NRPM by contradiction. Suppose
  $\exists \x \in \simplex^{|\node|}_\rho$ with
  $\x_i, \z^*_i(\x) \in [0,\theta_i]$ and $(i,j) \in \arc$ such that
  $\lev_i(\theta_i) \geq \lev_j(0)$ but $\d^*_{ij} > 0$ for some
  $\d^*$ that optimizes $\prob_3$. Then by
  \eqref{eq:level_order_nrpm}, we get
\begin{align*}
	\lev_i(\theta) \leq \lev_i(\z^*_i) \leq \lev_j(\z^*_j) < \lev_j(0) \,.
\end{align*}
The last strict inequality can be obtained by combining the feasibility constraints \eqref{eq:nrpm-flow-balance} for $j$, \eqref{eq:nrpm-xi-reallocation} and \eqref{eq:nrpm-non-neg} of $\prob_3$; along with the assumption that $\d^*_{ij} > 0$. But this is a contradiction and hence the initial assumption was incorrect. This completes the proof.
\end{proofa}

\subsubsection{Proof of Theorem \ref{th:upper_bound}} \label{proof:upper_bound}
\begin{proofa}
  From Remark \ref{rem:graph_red_init_st} it is clear that $\theta$ estimated using $\x^0$ in
  \texttt{reduceGraph}() has the property that $\x(t) \leq \theta$ and
  $\z^*(\x(t)) \leq \theta$, $\forall t \geq 0$ and hence the evolution on $\graph$ is same as that on $\redflowx$.

Now, the following arguments hold for SSD, NBRD and NRPM and we address them in one go. Suppose $\bx$ is the steady state state from the initial condition $\x^0$. Then notice that $\bx$ can be expressed as $\bx_i = \w_i = \sum_{j \in \inreach^i} \rvec_{ji}$, $\forall i \in \node$ for some $(\w,\rvec)$ which is a feasible solution of $\prob_4$. Thus $\umax \geq U(\bx)$. This completes the proof.
\end{proofa}

\subsubsection{Proof of Lemma \ref{lem:eventually_empty}} \label{proof:eventually_empty}
\begin{proofa}
\emph{(By contradiction)} Suppose $\exists \, i \in \rednodex$, such that $(i,j) \in \redarcx$ and $(j,i) \notin \redarcx$ but $\overline{\x}_i \neq 0$ for some $\bx \in \poslim$. Then,
\begin{equation*}
	\lev_j(\overline{\x}_j) \stackrel{\text{a}}{\geq} \lev_j(\theta_j) \stackrel{\text{b}}{\geq} \lev_i(0) \stackrel{\text{c}}{>} \lev_i(\overline{\x}_i) \stackrel{\text{d}}{\geq} \lev_j(\overline{\x}_j) \, .
\end{equation*}
Here, inequalities a and c comes from the strict decreasing nature of $\lev_i(.)$ and $\lev_j(.)$. Inequality c is strict as $\overline{\x}_i \neq 0$. Inequality b comes from the fact that $(j,i) \notin \redarcx$. Inequality d comes from the fact that $(i,j) \in \redarcx$ and $\overline{\x} \in \ne^{|\node|}_\rho$. This is a contradiction and the claim of the lemma is hence true.
\end{proofa}

\subsubsection{Proof of Lemma \ref{lem:super_node_bound_up}} \label{proof:super_node_bound_up}
\begin{proofa}
Consider an arbitrary but fixed super node $q \in \Lambda$. First we address the trivial cases. If $\M^q = \snode^q$, then \eqref{eq:snode_up_bound} is satisfied with $\rhom_q = 0$ by Lemma \ref{lem:eventually_empty}. Also as the total population is \emph{one}, then $\rhom_q = 1$ is always a correct upper bound on the total population fraction in $q$.

Next we address the case when $\snode^q \supset \M^q \neq \varnothing$. As $\hill^q$ is a DQCH component of $\redflowx$, by Lemma \ref{lem:same_lev} if $\bx_i \neq 0$ for some $i \in \snode^q \setminus \M^q$ then $\lev_i(\bx_i) \geq \lev_j(0) = \a_j$, $\forall j \in \M^q$ and hence $\lev_i(\bx_i) \geq \max_{j \in \M^q} \a_j = \a^\text{max}_q$. Thus, \eqref{eq:snode_up_bound} comes from applying this condition to all nodes $i \in \snode^q \setminus \M^q$ such that $\bx_i \neq 0$.
\end{proofa}

\subsubsection{Proof of Lemma \ref{lem:p1_concave_rho}} \label{proof:p1_concave_rho}
\begin{proofa}
\emph{By contradiction} Suppose $\exists \, \rho_1,\rho_2 \in \real_+$ and $\sigma \in (0,1)$ such that
\begin{equation}
	\optu_{\dumnode}\Big(\sigma\rho_1 + (1-\sigma)\rho_2 \Big) < \sigma\optu_{\dumnode}(\rho_1) + (1-\sigma)\optu_{\dumnode}(\rho_2)\,.
	\label{eq:last_contradict}
\end{equation}
We ignore the trivial cases where $\sigma = 0$ or $\sigma = 1$. Let $\x^1$ and $\x^2$ be the optimizers of $\prob_1(\dumnode,\rho_1)$ and $\prob_1(\dumnode,\rho_2)$ respectively. Thus, $\optu_{\dumnode}(\rho_1) = U(\x^1)$ and $\optu_{\dumnode}(\rho_2) = U(\x^2)$. Then it is easy to see that $\sigma \x^1 + (1-\sigma) \x^2$ is a feasible solution of $\prob_1(\dumnode,\sigma \rho_1 + (1-\sigma) \rho_2)$. Then,
\begin{equation*}
	\begin{split}
		& U \Big( \sigma \x^1 + (1-\sigma) \x^2 \Big) \leq \optu_{\dumnode} \Big( \sigma \rho_1 + (1-\sigma) \rho_2 \Big)\\
		& < \sigma\optu_{\dumnode}(\rho_1) + (1-\sigma)\optu_{\dumnode}(\rho_2) = \sigma U(\x^1) + (1-\sigma)U(\x^2) .
	\end{split}
\end{equation*}
Here the first inequality comes from the fact that $\optu_{\dumnode} \Big( \sigma \rho_1 + (1-\sigma) \rho_2 \Big)$ is the optimum of $\prob_1(\dumnode,\sigma \rho_1 + (1-\sigma) \rho_2)$. The second strict inequality comes from \eqref{eq:last_contradict}. This contradicts the strict concave nature of $U(.)$ and hence completes the proof. 

\end{proofa}

\subsubsection{Proof of Theorem \ref{th:lower_bound}} \label{proof:lower_bound}
\begin{proofa}
Note that the main steps of the process are illustrated in Algorithm \ref{algo:u_min}. By Remark \ref{rem:graph_red_init_st}, we know that $\redflowx$ in Step \ref{stp:final_algo_red_graph} computes an ICRG of $\graph$ and that the evolution of the population is unaffected by it. Next Step \ref{stp:final_algo_macp} computes a MAC partition of $\redflowx$, the accuracy of which is demonstrated in Remark \ref{rem:compute_macp}. $\Gamma$ in Step \ref{stp:final_algo_macsg} is the corresponding MAC-SG. By Lemma \ref{lem:super_node_bound_up}, $\rhom_q$ is an upper bound on the population fraction in each super node for every state in the positive limit set (say $\poslim$) of the trajectory starting from $\x^0$. Thus, every $\bx \in \poslim$ can be written as a feasible solution for $\prob_5$. Moreover, inside each super node, the steady state social utility is given by $\prob_1$. Hence the lower bound holds. 
\end{proofa}

\end{document}